\documentclass[10pt,a4paper]{article}
\usepackage{graphicx}
\usepackage{subcaption}
\usepackage{xcolor}
\usepackage{comment}
\usepackage[utf8]{inputenc}
\usepackage{comment}

\usepackage{hyperref}
\hypersetup{
    colorlinks=true,
    linkcolor=blue,
    citecolor = magenta,
    filecolor=magenta,      
    urlcolor=cyan,
}

\usepackage{siunitx,booktabs}
\usepackage{arydshln}
\usepackage{setspace,amsmath,amsthm,amsfonts,amssymb,mathtools,bm,algorithm,algpseudocode,
longtable,multirow,geometry,mathrsfs,tablefootnote,threeparttable,dsfont}

\usepackage[font=sf, labelfont={sf,bf}, margin=1cm]{caption}

\numberwithin{Theorem}{section}
\numberwithin{Definition}{section}
\numberwithin{Lemma}{section}
\numberwithin{Algorithm}{section}
\numberwithin{equation}{section}

\newtheorem{remark}{Remark}
\newcommand\scalemath[2]{\scalebox{#1}{\mbox{\ensuremath{\displaystyle #2}}}}

\makeatletter
\def\@cline#1-#2\@nil{%
  \omit
  \@multicnt#1%
  \advance\@multispan\m@ne
  \ifnum\@multicnt=\@ne\@firstofone{&\omit}\fi
  \@multicnt#2%
  \advance\@multicnt-#1%
  \advance\@multispan\@ne
  \leaders\hrule\@height\arrayrulewidth\hfill
  \cr
  \noalign{\nobreak\vskip-\arrayrulewidth}}
\makeatother
\begin{document}
	\title{An active-set method for sparse approximations\\
	Part II: General piecewise-linear terms}
\author{Spyridon Pougkakiotis \and Jacek Gondzio \and Dionysios S. Kalogerias}

\maketitle

			%\begin{changemargin}{0.8cm}{0.8cm} 
\begin{abstract}
\par In this paper we present an efficient active-set method for the solution of convex quadratic programming problems with general piecewise-linear terms in the objective, with applications to sparse approximations and risk-minimization. The method exploits the structure of the piecewise-linear terms appearing in the objective in order to significantly reduce its memory requirements, and thus improve its efficiency. We showcase the robustness of the proposed solver on a variety of problems arising in risk-averse portfolio selection, quantile regression, and binary classification via linear support vector machines. We provide computational evidence to demonstrate, on real-world datasets, the ability of the solver of efficiently handling a variety of problems, by comparing it against an efficient general-purpose interior point solver as well as a state-of-the-art alternating direction method of multipliers. This work complements the accompanying paper [``\emph{An active-set method for sparse approximations. Part I: Separable $\ell_1$ terms}", \emph{S. Pougkakiotis,  J. Gondzio, D. S. Kalogerias}], in which we discuss the case of separable $\ell_1$ terms, analyze the convergence, and propose general-purpose preconditioning strategies for the solution of its associated linear systems.
\end{abstract}
\section{Introduction}
\par In this paper we are interested in the solution of the following optimization problem
\begin{equation} \label{primal problem original} 
\begin{split}
\underset{x \in \mathbb{R}^n}{\text{min}}&\quad \left\{c^\top x + \frac{1}{2}x^\top Q x + \sum_{i = 1}^l \left( \left(C x + d\right)_{i}\right)_+ + \|Dx\|_1 + \delta_{\mathcal{K}}(x)\right\}, \\ 
\textnormal{s.t.} &\quad \ A x = b,
\end{split}
\end{equation}
\noindent where $c \in \mathbb{R}^n$, $Q \in \mathbb{R}^{n \times n}$ is a positive semi-definite matrix, $C \in \mathbb{R}^{l\times n}$, $d \in \mathbb{R}^l$, $A \in \mathbb{R}^{m\times n}$ is a linear constraint matrix with $b\in \mathbb{R}^m$ a given right-hand side, $D \in \mathbb{R}^{n\times n}$ is a diagonal positive semi-definite (``weight") matrix, and $\mathcal{K}$ is a closed convex set $\mathcal{K} \coloneqq \{x \in \mathbb{R}^n \colon x \in \left[a_l,a_u\right]\}$ with $a_l,\ a_u \in \mathbb{R}^n$, such that $(a_l)_i \in \mathbb{R}\cup\{-\infty\},\ (a_u)_i \in \mathbb{R}\cup \{+\infty\}$ for all $i = 1,\ldots,n$. Additionally, $(\cdot)_+ \equiv \max\{\cdot,0\}$, while $\delta_{\mathcal{K}}(\cdot)$ is an indicator function for the set $\mathcal{K}$, that is, $\delta_{\mathcal{K}}(x) = 0$ if $x \in \mathcal{K}$ and $\delta_{\mathcal{K}}(x) = \infty$, otherwise.
\par We introduce an auxiliary variable $w \in \mathbb{R}^l$, and reformulate \eqref{primal problem original} in the following form:
\begin{equation} \label{primal problem} \tag{P}
\begin{split}
\underset{(x,w)\ \in\ \mathbb{R}^n\times\mathbb{R}^l}{\text{min}} & \left\{c^\top   x + \frac{1}{2}x^\top Q x +  \sum_{i = 1}^l \left(w_i\right)_+ + \|Dx\|_1 + \delta_{\mathcal{K}}(x)\right\},\\
\textnormal{s.t.}\quad\ \ &\ Cx + d - w = 0_l,\\
&\ Ax = b.
\end{split}
\end{equation}

\begin{remark} \label{remark: l1 norm reformulation}
\par Let us notice that the model in \eqref{primal problem} can readily accommodate terms of the form of $\|Cx + d\|_1$ where $C \in \mathbb{R}^{l\times n}$, and $d \in \mathbb{R}^l$. Indeed, letting $c = - C^\top \mathds{1}_l$ and adding the constant term $-\mathds{1}_l^\top d$ in the objective of \eqref{primal problem}, we notice that
\[\|Cx + d\|_1 \equiv -\mathds{1}_l^\top (Cx + d) + \sum_{j=1}^l \left(2(Cx + d)_j\right)_+, \]
where $\mathds{1}_l \coloneqq (1,\ldots,1)^\top \in \mathbb{R}^l$. Similarly, any piecewise-linear term of the form 
\[\sum_{i=1}^l\max\left\{\left(C_1 x+ d_1\right)_i, \left(C_2 x+ d_2\right)_i\right\},\]
\noindent where $C_1,\ C_2 \in \mathbb{R}^{l\times n}$ and $d_1,\ d_2 \in \mathbb{R}^l$, can also be readily modeled. Indeed, setting $c = C_2^\top\mathds{1}$ and adding the term $d_2^\top \mathds{1}_l$ in the objective yields
\[ \mathds{1}_l^\top \left(C_2 x + d_2\right) + \sum_{i=1}^l \left(\left(C_1 x + d_1 - C_2x - d_2\right)_i \right)_+ \equiv \sum_{i = 1}^l \max\left\{\left(C_1 x + d_1\right)_i,\left(C_2 x + d_2\right)_i \right\}.\]
\noindent Finally, it is important to note that model \eqref{primal problem} allows for multiple piecewise-linear terms of the form $\max\{Cx + d,0_l\}$, $\|Cx + d\|_1$ or $\max\left\{C_1 x+ d_1, C_2 x+ d_2\right\}$, since we can always adjust $l$ to account for more than one terms. Hence, one can observe that \eqref{primal problem} is quite general and can be used to model a plethora of very important problems that arise in practice. 
\end{remark}
\par In light of the discussion in Remark \ref{remark: l1 norm reformulation}, it is easily observed that problem \eqref{primal problem} can model a plethora of very important problems arising in several application domains spanning, among others, operational research, machine learning, data science, and engineering. More specifically, various lasso and fussed lasso instances (with applications to sparse approximations for classification and regression \cite{SIREV:DeSimone_etal,JRSS:Zou}, portfolio allocation \cite{JBanFin:AlexColYi}, or medical diagnosis \cite{IPRN:Gramfort_etal}, among many others) can be readily modeled by \eqref{primal problem}. Additionally, various risk-minimization problems with linear random cost functions can be modeled by \eqref{primal problem} (e.g. see \cite{MathFin:LiNg,JRisk:RockUry,Prod:Silva_etal}). Furthermore, even risk-minimization problems with nonlinear random cost functions, which are typically solved via Gauss-Newton schemes (e.g. see \cite{MathProg:Burke}), often require the solution of sub-problems like \eqref{primal problem}. Finally, continuous relaxations of integer programming problems with applications to operational research (e.g. \cite{MathProg:Laz}) often take the form of \eqref{primal problem}. Given the multitude of problems requiring easy access to (usually accurate) solutions of \eqref{primal problem}, the derivation of efficient, robust, and scalable solution methods is of paramount importance. 
\par Problem \eqref{primal problem} can be solved by various first- or second-order methods. In particular, using a standard reformulation, by introducing several auxiliary variables, \eqref{primal problem} can be written as a convex quadratic programming (QP) one and efficiently solved by, among others, an \emph{interior point method} (IPM; e.g. \cite{SIREV:DeSimone_etal,COAP:PougkGond}), an \emph{alternating direction method of multipliers} (ADMM; e.g. \cite{SciComp:DengYin}), or a \emph{proximal point method} (e.g. \cite{COAP:KanzowLechner,SIAMJCO:Rock}). However, the reformulation of \eqref{primal problem} into a convex QP is not expected to lead to scalable solution methods, since the dimension of the problem is significantly increased, and hence an already large-scale instance might be very difficult to solve in this way. 
\par Alternatively, ADMM (or general splitting) schemes can be developed for the solution of \eqref{primal problem} without the need of additional auxiliary variables (see Section \ref{sec: warm start}). However, no first-order scheme would be able to consistently yield sufficiently accurate solutions (i.e. of 4-, 5- or 6-digit accuracy). If such a solution is sought, we have to employ a semismooth Newton method (SSN; e.g. \cite{NLAA:ChenQi,MathOR:Han,CAM:MartQi}), or a combination of a proximal point method with an SSN scheme utilized for the solution of its associated sub-problems. SSN methods have been successfully applied to a plethora of problems, however, their success is heavily reliant on the properties of the problem at hand (e.g. the rank of the linear constraints, or the conditioning of the Hessian). On the other hand, the combination of the proximal point method with the SSN can circumvent these issues, since the associated nonsmooth sub-problems can be guaranteed to be well-defined and well-conditioned.
\par Various such solvers have been developed and analyzed in the literature. For example, the authors in \cite{SIAMOpt:ZhaoSunToh} developed a dual augmented Lagragian scheme combined with an SSN method for the solution of semidefinite programming problems, and obtained very promising results. This scheme was then utilized for the solution of linear programming problems in \cite{SIAMOpt:Lietal}, and for lasso-regularized problems in \cite{SIAMOpt:Lietal2}. A similar primal-dual approach for $\ell_1$-regularized convex quadratic programming problems was developed and analyzed in our accompanying paper \cite{arxiv:PougkGondzio} and was shown to be especially efficient for the solution of elastic-net linear regression and $L^1$-regularized partial differential equation constrained optimization problems. In fact, the proposed active set method developed in this work is a direct extension of the method given in \cite{arxiv:PougkGondzio}, altered in a specific way so that it can efficiently handle most piecewise-linear terms that appear in practice, via restricting its memory requirements.
\par Indeed, we showcase that each of the nonsmooth terms in the objective of \eqref{primal problem} can be utilized for reducing the memory requirements of the proposed method. In particular, when computing the \emph{Clarke subdifferential} (\cite{JWS:Clarke}) of an appropriate augmented Lagrangian penalty function, we can show that the Clarke derivatives of such piecewise-linear terms can act as projectors. As a result, we obtain an active-set scheme that reduces the sub-problems' dimensions significantly, leading to better performance and reduced memory requirements. In particular, we observe that a thresholding operation (originating from the presence of $\ell_1$ terms in the objective) determines which of the variables $x$ are inactive, allowing us to throw away entire columns from matrices $A$ and $C$ when solving the associated sub-problems. Furthermore, the $\max\{\cdot,0\}$ terms in the objective determine which of the rows of $C$ are non-important, allowing us to further eliminate such rows. 
\par We showcase the robustness and the efficiency of the resulting active-set scheme on various optimization problems arising in risk-averse portfolio selection, quantile regression, and binary classification via linear support vector machines. In each of these cases the proposed scheme is compared against the robust polynomially convergent regularized interior point method developed and analyzed in \cite{COAP:PougkGond}, as well as the well-known ADMM-based OSQP solver (\cite{osqp}). We demonstrate the reduced memory requirements of the active set scheme (and hence its improved scalability), as compared to interior point and ADMM alternatives applied to QP reformulations, while showcasing its efficiency and robustness. 
\paragraph{Structure of the paper} In Section \ref{sec: active-set method} we briefly present the proposed inner-outer active-set method. In particular, in Section \ref{subsec: AL penalties} we derive a proximal method of multipliers (outer scheme) for the solution of \eqref{primal problem}, assuming that we can find an $\epsilon$-optimal solution of the associated sub-problems. Then, in Section \ref{subsec: SSN method}, we briefly present a semismooth Newton method (inner scheme) for finding approximate solutions to sub-problems arising from the proximal method of multipliers. Focusing on the structure of problem \eqref{primal problem}, by selecting appropriate Clarke derivatives, we show that the proposed inner-outer method is in fact an active-set scheme, the associated linear systems of which are well-conditioned and stable. In Section \ref{subsection: extension to robust optimization}, we discuss an extension of the method for dealing with problems having arbitrary piecewise-linear terms in the objective (with applications to, among others, robust optimization) that are not currently covered within the paper.
\par Subsequently, in Section \ref{sec: warm start}, we derive a proximal alternating direction method of multipliers for the approximate solution of \eqref{primal problem} in order to obtain good initial estimates for the primal and dual variables of the problem. This can then be used to warm-start the proposed second-order algorithm. A good starting point for the algorithm could mean that only a small portion of the problem data matrices is used at each inner-outer iteration of the scheme, while the outer method is expected to achieve its local linear (and potentially superlinear) convergence rate in a small number of inner-outer iterations. 
\par In Section \ref{sec: Applications}, we showcase the efficiency and robustness of the approach on several important real-life applications arising in risk-averse portfolio optimization, statistical regression, and binary classification via linear support vector machines. Finally, we discuss our conclusions in Section \ref{sec: Conclusions}.
\paragraph{Notation} Given a vector $x$ in $\mathbb{R}^n$, $\|x\|$ denotes its Euclidean norm. Given a closed set $\mathcal{K} \subset \mathbb{R}^n$, $\Pi_{\mathcal{K}}(x)$ denotes the Euclidean projection onto $\mathcal{K}$, that is $\Pi_{\mathcal{K}}(x) \coloneqq \arg\min\{\|x-z\| \colon z \in \mathcal{K}\}$. Given a closed set $\mathcal{K}$, we write $\textnormal{dist}(z,D) \coloneqq \inf_{z'\in \mathcal{K}}\|z-z'\|$. Given a convex function $p \colon \mathbb{R}^n \mapsto \mathbb{R}$, we define the proximity operator as $\textbf{prox}_p (u) \coloneqq \arg\min_x \left\{ p(x) + \frac{1}{2}\|u-x\|^2\right\}$. Given an index set  $\mathcal{D}$, $|\mathcal{D}|$ denotes its cardinality. Given a rectangular matrix $A \in \mathbb{R}^{m \times n}$ and an index set $\mathcal{B} \subseteq \{1,\ldots,n\}$, we denote the columns of $A$, the indices of which belong to $\mathcal{B}$, as $A_{\mathcal{B}}$. Given a matrix $A \in \mathbb{R}^{m\times n}$, and two index sets $\mathcal{B}_1 \subseteq \{1,\ldots,m\}$ and $\mathcal{B}_2 \subseteq \{1,\ldots,n\}$, we denote the subset of rows and columns of $A$, the indices of which belong to $\mathcal{B}_1,\ \mathcal{B}_2$ respectively, as $A_{\left(\mathcal{B}_1,\mathcal{B}_2\right)}$. Finally, given an arbitrary vector $d$ with $n$ components as well as some indexes $1\leq i_1 \leq i_2 \leq n$, we denote by $d_{i_1:i_2}$ the vector containing the $i_1$-st up to the $i_2$-nd component of this vector. To avoid confusion, indexes $i$ are always used to denote entries of vectors or matrices, while $k$ and $j$ are reserved to denote iteration counters (outer and inner, respectively).

\section{An active-set method} \label{sec: active-set method}
\par In this section we derive an active-set method for the solution of \eqref{primal problem}. The algorithm is an inner-outer scheme, which results by combining an outer proximal method of multipliers, and an inner semismooth Newton method for the solution of the PMM sub-problems. Following the discussion in \cite[Section 2]{arxiv:PougkGondzio}, we briefly derive the outer PMM scheme. Then, we briefly present the inner semismooth Newton method and discuss the solution of its associated linear systems, which is the main bottleneck of the algorithm.
\subsection{Outer scheme: A primal-dual proximal method of multipliers} \label{subsec: AL penalties}
\par In this section, following and extending the developments in \cite{arxiv:PougkGondzio}, we derive a primal-dual proximal method of multipliers for the approximate solution of \eqref{primal problem}. A convergence analysis is not provided, and the reader is referred to \cite[Section 2]{arxiv:PougkGondzio}, for an outline of such an analysis (since it applies directly to the case under consideration).
\par Given a penalty parameter $\beta > 0$, and some dual multiplier estimates $y,\ z$, we follow \cite[Section 2]{arxiv:PougkGondzio} to obtain the augmented Lagrangian corresponding to \eqref{primal problem}, which reads
\begin{equation} \label{final augmented lagrangian of the primal}
\begin{split}
\mathcal{L}_{\beta}(x, w; y, z) \coloneqq\ &\ c^\top  x + \frac{1}{2}x^\top Q x + \sum_{i = 1}^l \left(w_i\right)_+ + \|Dx\|_1 - y^\top \left(\begin{bmatrix}Cx+d- w  \\ Ax-b \end{bmatrix}\right) \\ &\quad \ + \frac{\beta}{2}\left\|\begin{bmatrix}Cx + d- w\\ Ax-b\end{bmatrix}\right\|^2 - \frac{1}{2\beta}\|z\|^2  + \frac{1}{2\beta}  \|z + \beta x - \beta \Pi_{\mathcal{K}}(\beta^{-1}z + x)\|^2.
\end{split}
\end{equation}
\noindent Indeed, this can be verified utilizing certain standard properties of Fenchel duality and of the proximity operator (see \cite[Section 2]{arxiv:PougkGondzio}). 
\par During iteration $k\geq 0 $ of the proximal method of multipliers, we have the estimates $(x_k,y_k,z_k)$ as well as the penalty parameters $\beta_k,\ \rho_k$, such that $\rho_k \coloneqq \frac{\beta_k}{\tau_k}$, where $\{\tau_k\}_{k=0}^{\infty}$ is a non-increasing positive sequence. For simplicity of exposition, let $g\left(x,w\right) \coloneqq g_1(x) + g_2(w)$, where $g_1(x) \coloneqq \|Dx\|_1,$ and $g_2(w) \coloneqq  \sum_{i = 1}^l \left(w_i\right)_+.$ We consider the following regularized continuously differentiable function:
\begin{equation*}
\begin{split} 
\scalemath{1}{\phi(x,w) \equiv\  \phi_{\beta_k,\rho_k}(x,w;x_k,y_k,z_k) \coloneqq }&\ \scalemath{1}{\mathcal{L}_{\beta_k}(x,w;y_k,z_k) + \frac{1}{2\rho_k}\left\|x-x_k\right\|^2  - g\left(x,w\right)}. 
\end{split}
\end{equation*}
\noindent Notice that we introduce primal proximal regularizer only for variable $x$, and not for $w$. This is a very important algorithmic choice that departs from the developments in \cite{arxiv:PougkGondzio}. We want to develop a memory-efficient active-set method, and for that reason we avoid introducing a proximal term for the auxiliary variables $w$. This choice, which does not hinder the stability of the proposed approach, leads to better memory efficiency of the resulting active-set scheme, by promoting a sparsification of the associated linear systems solved at each inner-outer iteration (this point will become clear in Section \ref{subsubsec: SSN linear systems}). 
\par The minimization of the proximal augmented Lagrangian function can be written as 
\[\underset{x,w}{\min}\ \left\{\mathcal{L}_{\beta_k}(x,w;y_k,z_k) + \frac{1}{2\rho_k}\left\|x-x_k\right\|^2\right\} \equiv  \min_{x,w} \left\{\phi(x,w)+ g(x,w)\right\},\]
\noindent and thus we need to find $(x^*,w^*) \in \mathbb{R}^n \times \mathbb{R}^l$ such that 
\[\scalemath{1}{\left(\nabla \phi(x^*,w^*)\right)^\top \left((x,w) - (x^*,w^*)\right) + g(x,w) - g(x^*,w^*) \geq 0,\quad \forall\ (x,w)\in \mathbb{R}^n\times  \mathbb{R}^l}.\]
\noindent To that end, we observe that
\begin{equation*}
\begin{split}
\nabla_x \phi(x,t,w) =&\  c + Qx - \begin{bmatrix}C^\top & A^\top\end{bmatrix} y_k + \beta_k \begin{bmatrix}C^\top & A^\top\end{bmatrix}  \left(\begin{bmatrix}Cx +d-  w\\ Ax-b\end{bmatrix}\right) \\ & \qquad + (z_k + \beta_k x) - \beta_k \Pi_{\mathcal{K}}(\beta_k^{-1}z_k + x) + \rho_k^{-1} (x-x_k), \\
\nabla_w \phi(x,t,w) =&\ 
 (y_k)_{1:l} - \beta_k \left( Cx + d - w\right).
\end{split}
\end{equation*}
\noindent By introducing the dual variable $y \coloneqq y_k - \beta_k \left(\begin{bmatrix}Cx +d- w \\ Ax - b\end{bmatrix}\right) \in \mathbb{R}^{l+m}$, the optimality conditions of $\underset{x,w}{\min}\ \left\{\phi(x,w)+ g(x,w)\right\}$ can be written as
\begin{equation} \label{PD-SSN-PMM Proximal Optimality Conditions}
(0_{n+l},0_{l+m}) \in  \mathcal{M}_{\beta_k,\rho_k}(x,w,y; x_k,y_k,z_k) \equiv\mathcal{M}_{k}(x,w,y),
\end{equation}
\noindent where 
\begin{equation*}
\begin{split}
 \mathcal{M}_{k}(x,w,y) \coloneqq \left\{\scalemath{1}{(u',v') \colon  u' \in r_{k}(x,y) + \begin{bmatrix} \partial g_1(x)\\ \partial g_2(w) \end{bmatrix},\  v' = \begin{bmatrix}Cx + d-  w\\ Ax-b\end{bmatrix} + \beta_k^{-1}(y-y_k)}\right\},
 \end{split}
 \end{equation*}
\[r_{k}(x,y) \coloneqq  \begin{bmatrix} \scalemath{1}{c + Qx - \begin{bmatrix} C^\top & A^\top \end{bmatrix} y  + (z_k + \beta_k x) - \beta_k \Pi_{\mathcal{K}}(\beta_k^{-1}z_k + x)} + \rho_k^{-1}(x-x_k) \\
 y_{1:l}
\end{bmatrix}.\]
\noindent We now describe the proposed proximal method of multipliers in Algorithm \ref{PMM algorithm}.
\renewcommand{\thealgorithm}{PMM}

\begin{algorithm}[!ht]
\caption{(\emph{proximal method of multipliers})}
    \label{PMM algorithm}
    \textbf{Input:}  $(x_0,w_0,y_0,z_0) \in \mathbb{R}^n \times \mathbb{R}^l \times \mathbb{R}^{l+m} \times \mathbb{R}^n$, $\beta_0,\ \beta_{\infty} > 0$, $\{\tau_k\}_{k=0}^{\infty}$, such that $\tau_k \geq \tau_{\infty} > 0$, $\forall\ k \geq 0$.
\begin{algorithmic}
\State Choose a sequence of positive numbers $\{\epsilon_k\}$ such that $\epsilon_k \rightarrow 0$. 
\For {($k = 0,1,2,\ldots$)}
\State Find $(x_{k+1},w_{k+1},y_{k+1})$ such that:
\begin{equation} \label{primal-dual PMM main sub-problem}
\varepsilon_k \coloneqq  \textnormal{dist}\left(0_{n+2l+m},\mathcal{M}_{k}\left(x_{k+1},w_{k+1},y_{k+1}\right)\right)  \leq \epsilon_k,
\end{equation}
\State where letting $\hat{r} = r_{k}(x_{k+1},y_{k+1})$ and $\mathcal{U} = \left\{u \in \mathbb{R}^{n+l} \colon u \in \partial g(x_{k+1},w_{k+1})\right\},$ we have
\[\scalemath{1}{\varepsilon_k = \left\|\begin{bmatrix}\hat{r}+ \Pi_{\mathcal{U}}\left(-\hat{r}\right)\\
\begin{bmatrix} C \\ A \end{bmatrix} x_{k+1} + \begin{bmatrix} d \\ - b\end{bmatrix} -\begin{bmatrix} I_l \\ 0_{m,l}\end{bmatrix}w_{k+1} + \beta_k^{-1} (y_{k+1} - y_k)
\end{bmatrix}\right\|.} \]
   \begin{flalign}  \label{primal-dual PMM z-update}
\ \ \ \ z_{k+1} &= \  (z_k + \beta_k x_{k+1}) - \beta_k\Pi_{\mathcal{K}}\big(\beta_k^{-1}z_k + x_{k+1} \big).&&
\end{flalign} 
\begin{flalign} \ \ \ \ \beta_{k+1} &\nearrow \beta_{\infty} \leq \infty, \quad \rho_{k+1} = \frac{\beta_{k+1}}{\tau_{k+1}}.&&
\end{flalign}
\EndFor
\State \Return $(x_{k+1},w_{k+1},y_{k+1},z_{k+1})$.
\end{algorithmic}
\end{algorithm}
\par Let us notice that given a triple $(\tilde{x},\tilde{w},\tilde{y})$, we can easily evaluate the distance in \eqref{primal-dual PMM main sub-problem}, due to the piecewise linear structure of the associated function $g(\cdot)$. A trivial extension of the analysis in \cite[Section 3]{arxiv:PougkGondzio} yields that Algorithm \ref{PMM algorithm} is globally convergent under the standard assumption of primal-dual feasibility. Furthermore, given some additional conditions on the error sequence $\{\epsilon_k\}$ one can show that in fact Algorithm \ref{PMM algorithm} achieves a local linear convergence rate (which becomes superlinear if $\beta_{k} \rightarrow \infty$ at a suitable rate). Finally, the algorithm exhibits a global linear convergence rate, assuming that the starting point is chosen properly. For more details, we refer the reader to \cite[Theorems 2.2, 2.3]{arxiv:PougkGondzio}.

\subsection{Inner scheme: A semismooth Newton method} \label{subsec: SSN method}
\par Next, we employ a semismooth Newton (SSN) method to solve problem \eqref{primal-dual PMM main sub-problem} appearing in Algorithm \ref{PMM algorithm}, and verify that the resulting inner-outer scheme admits an active-set interpretation.  Given the estimates $(x_k,y_k,z_k)$ as well as the penalties $\beta_k,\ \rho_k$, we apply SSN to approximately solve \eqref{PD-SSN-PMM Proximal Optimality Conditions}. Given any bounded positive penalty $\zeta_k > 0$, the optimality conditions in \eqref{PD-SSN-PMM Proximal Optimality Conditions} can be written as
\begin{equation} \label{PD-SSN-PMM Proximal Optimality Conditions-SMOOTH}
\widehat{\mathcal{M}}_{k}\left(x,w,y\right) \coloneqq \zeta_k\begin{bmatrix}
\zeta_k^{-1}\left(\begin{bmatrix} x \\ w \end{bmatrix} - \textbf{prox}_{\zeta_k g}\left(x- \zeta_k r_{k}(x,y) \right) \right)
\\ 
\begin{bmatrix} C \\ A \end{bmatrix} x + \begin{bmatrix} d \\ - b\end{bmatrix} -\begin{bmatrix} I_l \\ 0_{m,l}\end{bmatrix}w + \beta_k^{-1} (y - y_k)
\end{bmatrix}  = \begin{bmatrix}
0_{n+l}\\
0_{l+m}
\end{bmatrix}.
\end{equation}
\noindent We set $x_{k_0} = x_k$, $y_{k_0} = y_k$, and at every iteration $k_j$ of SSN, we solve 
\begin{equation} \label{Primal-dual SSN system}
\underbrace{\scalemath{1}{\begin{bmatrix}
H_{k_j}   & 0_{n,l} & -\zeta_k B_{g_1,k_j}\begin{bmatrix} C^\top & A^\top \end{bmatrix}\\
0_{l,n} &  \left(I_l - B_{g_2,k_j}\right) &  \zeta_k \begin{bmatrix}B_{g_2,k_j} & 0_{l,m}\end{bmatrix}\\
\zeta_k\begin{bmatrix} C\\ A\end{bmatrix} &  \zeta_k\begin{bmatrix} -I_l \\ 0_{m,l}\end{bmatrix} & \zeta_k\beta_k^{-1} I_{l+m}
\end{bmatrix}}}_{{M_{k_j}}}\begin{bmatrix}
d_x\\
d_w\\
d_y
\end{bmatrix}= - \widehat{\mathcal{M}}_{k}\left(x_{k_j},w_{k_j},y_{k_j}\right) ,
\end{equation}
\noindent where we have introduced the notation 
\[H_{k_j} \coloneqq I_n - B_{g_1,k_j} + \zeta_k\beta_k B_{g_1,k_j}\left(\left(1+\rho_k^{-1}\beta_k^{-1}\right)I_n - B_{\delta,k_j} +\beta_k^{-1}Q \right),\] assuming that we are given some arbitrary matrices
\begin{equation} \label{eqn: Clarke subdifferential matrices}
\begin{split}
& B_{\delta,k_j} \in \partial^C_x \Pi_{\mathcal{K}}\left(\beta_k^{-1}z_k + x_{k_j}\right),\\
& B_{g_1,k_j} \in \partial^C_x \left(\textbf{prox}_{\zeta_k g_1}\left(x_{k_j}- \zeta_k(r_{k}(x_{k_j},y_{k_j}))_{1:n}\right)\right),\\
& B_{g_2,k_j} \in \partial^C_w \left(\textbf{prox}_{\zeta_k g_2}\left(w_{k_j} - \zeta_k  (y_{k_j})_{1:l}\right) \right).
\end{split}
\end{equation}
\noindent The symbol $\partial_x^C(\cdot)$ denotes the Clarke subdifferential of a function (see \cite{JWS:Clarke}). In our case, any element of the Clarke subdifferential is a \emph{Newton derivative} (see \cite[Chapter 13]{arXiv:ClasValk}), since $r_{k}(\cdot,\cdot)$ and $g(\cdot,\cdot) \equiv g_1(\cdot) + g_2(\cdot)$ are \emph{piecewise continuously differentiable} and \emph{regular functions}. For any $u \in \mathbb{R}^n$ and any $i \in \{1,\ldots,n\}$, it holds that
\begin{equation*}
\partial_{u_i}^C \left(\Pi_{[a_i,b_i]}(u_i)\right)= \begin{cases} 
\{1\}, &\qquad \textnormal{if}\quad u_i \in (a_i,b_i),\\
\{0\}, &\qquad \textnormal{if}\quad u_i \notin [a_i,b_i],\\
[0,1], &\qquad \textnormal{if}\quad u_i \in \{a_i,b_i\}.
\end{cases}
\end{equation*}
\noindent Since $g_1(x) = \|Dx\|_1$ and $D\succeq 0_n$ is diagonal, we have that for all $u \in \mathbb{R}^n$ and any $i \in \{1,\ldots,n\}$,
\begin{equation*}
\left( \textbf{prox}_{\zeta_k g_1}\left(u\right)\right)_i = \textbf{soft}_{\left[-\zeta_kD_{(i,i)},\zeta_k D_{(i,i)}\right]}\left(u_i\right) \equiv \max\left\{\lvert u_i \rvert - \zeta_k D_{(i,i)},0 \right\}\textnormal{sign}(u_i).
\end{equation*}
\noindent Finally, since $g_2(w) = \max\{w,0\}$, we have that for all $u \in \mathbb{R}^l$ and any $i = 1,\ldots,l$,
\begin{equation*}
\left( \textbf{prox}_{\zeta_k g_2}\left(u\right)\right)_i = \textbf{soft}_{\left[0,\zeta_k\right]}\left(u_i\right) \equiv \max\{u_i-\zeta_k,0\} +\min\{u_i,0\}.
\end{equation*}
\noindent Then, we can show (e.g. see \cite[Example 14.9]{arXiv:ClasValk}) that
\begin{equation*}
\scalemath{0.95}{\left(\partial_{u}^C \left(\textbf{prox}_{\zeta_k g_1}\left(u\right)\right)\right)_i= \begin{cases} 
\{1\}, &\qquad \textnormal{if}\quad  \lvert u_i \rvert > \zeta_k D_{(i,i)} \ \textnormal{or}\ D_{(i,i)} = 0,\\
\{0\}, &\qquad \textnormal{if}\quad  \lvert u_i \rvert < \zeta_k D_{(i,i)},\\
[0,1], &\qquad \textnormal{if}\quad  u_i  = \zeta_k D_{(i,i)}, 
\end{cases}\  \textnormal{for all } i \in \{1,\ldots,n\},}
\end{equation*}
\noindent and
\begin{equation*}
\left(\partial_{u}^C \left(\textbf{prox}_{\zeta_k g_2}\left(u\right)\right)\right)_i= \begin{cases} 
\{1\}, &\qquad \textnormal{if}\quad  u_i > \zeta_k\ \textnormal{or}\ u_i < 0,\\
\{0\}, &\qquad \textnormal{if}\quad 0 < u_i < \zeta_k,\\
[0,1], &\qquad \textnormal{if}\quad  u_i  = \zeta_k\ \textnormal{or}\ u_i = 0,
\end{cases} \  \textnormal{for all } i \in \{1,\ldots,l\}.
\end{equation*}

\par We complete the derivation of the SSN, by defining a primal-dual merit function for globalizing the semismooth Newton method via a backtracking line-search scheme. Following \cite{arxiv:PougkGondzio}, we employ the following merit function
\begin{equation} \label{merit function for SSN globalization}
\Theta_k(x,w,y) \coloneqq  \left\|\widehat{M}_{k}\left(x,w,y\right)\right\|^{2}.
\end{equation}
\noindent This function is often used for globalizing SSN schemes applied to nonsmooth equations of the form of \eqref{PD-SSN-PMM Proximal Optimality Conditions-SMOOTH}. In Algorithm \ref{primal-dual SNM algorithm}, we outline a locally superlinearly convergent semismooth Newton method for the approximate solution of \eqref{primal-dual PMM main sub-problem}. The associated linear systems can be solved approximately (e.g. by Krylov subspace methods as in \cite{arxiv:PougkGondzio}), although in this work a suitable factorization scheme is utilized. 
\renewcommand{\thealgorithm}{SSN}
\begin{algorithm}[!ht]
\caption{(\emph{semismooth Newton method})}
    \label{primal-dual SNM algorithm}
    \textbf{Input:}  Let $\epsilon_k > 0$, $\mu \in (0,\frac{1}{2})$, $\delta \in (0,1)$, $x_{k_0} (= x_k)$, $w_{k_0} (= w_k)$,  $y_{k_0} (= y_k)$.
\begin{algorithmic}
\For {($j = 0,1,2,\ldots$)}
\State Choose some $M_{k_j}$, as in \eqref{Primal-dual SSN system} and solve
\[M_{k_j} \begin{bsmallmatrix}d_x\\d_w\\d_y \end{bsmallmatrix} = - \widehat{M}_{k}\left(x_{k_j}, w_{k_j},y_{k_j}\right).\]
\State (\emph{Line Search}) Set $\alpha_j = \delta^{m_j}$, where $m_j$ is the first non-negative integer for which:
\begin{equation*}
\Theta_k\left(\bar{x},\bar{w},\bar{y}\right) \leq \left(1 - 2 \mu \delta^{m_j}\right) \Theta_k\left(x_{k_j},w_{k_j},y_{k_j}\right)
\end{equation*}
\State \noindent where $\left(\bar{x},\bar{w},\bar{y}\right) = \left(x_{k_j} + \delta^{m_j}d_x, w_{k_j} + \delta^{m_j}d_w, y_{k_j} + \delta^{m_j} d_y\right)$.
\State Set $\left(x_{k_{j+1}},w_{k_{j+1}},y_{k_{j+1}}\right) = \left(\bar{x},\bar{w},\bar{y}\right)$.
\If {$\left(\textnormal{dist}\left(0_{n+2l+m},M_{k}\left(x_{k_{j+1}},w_{k_{j+1}},y_{k_{j+1}}\right)\right) \leq \epsilon_k\right)$}
\State \Return $(x_{k_{j+1}},w_{k_{j+1}},y_{k_{j+1}})$.
\EndIf
\EndFor
\end{algorithmic}
\end{algorithm}
\par Local superlinear convergence of Algorithm \ref{primal-dual SNM algorithm} follows directly from \cite[Thereom 3]{CAM:MartQi}, since the equation in \eqref{PD-SSN-PMM Proximal Optimality Conditions-SMOOTH} is \emph{BD-regular} (i.e. the Bouligand subdifferential at the optimum contains nonsingular matrices) by construction, as it models ``regularized" sub-problems arising from the proximal method of multipliers. On the other hand, if we assume that the directional derivative of \eqref{merit function for SSN globalization} is continuous at the optimum, then we can mirror the analysis in \cite{InverseProbs:HansRaasch} to show that Algorithm \ref{primal-dual SNM algorithm} is globally convergent. There is a wide literature on similar semismooth Newton schemes, and we refer the reader to \cite{NLAA:ChenQi,InverseProbs:HansRaasch,CAM:MartQi,MathOR:Qi} and the references therein for more details.
%\par The are various other alternative semismooth Newton schemes in the literature, some of which might exhibit better theoretical guarantees. For example, there have been developed approaches based on trust-region strategies (e.g. see \cite{JOTA:Akbari_etal,SIAMOpt:Christofetal,MathProg:Dennis_etal,OE:MannelRund}), or line-search strategies based on smooth penalty functions (such as the forward-backward envelope (FBE) \cite{SIAMOpt:Themelis_etal}). Following \cite{arxiv:PougkGondzio}, we have chosen to employ \eqref{merit function for SSN globalization} since it provides an active-set scheme that appears to be especially efficient and scalable in practice, while at the same time exhibiting robustness for a wide range of practical applications.
\subsubsection{The SSN linear systems} \label{subsubsec: SSN linear systems}
\par The major bottleneck of the previously presented inner-outer scheme is the solution of the associated linear systems given in \eqref{Primal-dual SSN system}. One could alter Algorithm \ref{primal-dual SNM algorithm} so that it does not require an exact solution. In turn this would allow for the utilization of preconditioned Krylov subspace solvers for the efficient solution of such systems (e.g. as in \cite[Section 3]{arxiv:PougkGondzio}). In particular, any preconditioner derived in \cite{arXiv:GondPougkPears} can be utilized for the proposed solver. However, for simplicity of exposition, we employ a standard factorization approach. As will become evident in Section \ref{sec: Applications}, the active-set nature of the proposed algorithm enables the use of factorization even for large-scale problems of interest, since most of the data are inactive at each inner-outer iteration $k_j$. Indeed, in the presence of multiple piecewise-linear terms in the objective, and assuming that the method is appropriately warm-started, one can ensure that most of the problem data are not incorporated when forming the active Hessian at each inner-outer iteration.
\par In what follows we derive the associated linear systems. Let $k \geq 0$ be an arbitrary iteration of Algorithm \ref{PMM algorithm}, and $j \geq 0$ an arbitrary iteration of Algorithm \ref{primal-dual SNM algorithm}. Firstly, let us notice that any element $B_{\delta} \in \partial_x^C\left(\Pi_{\mathcal{K}}\left(\cdot \right)\right)$ yields a Newton derivative (see \cite[Theorem 14.8]{arXiv:ClasValk}). The same applies for any element $B_{g_1} \in \partial_x^C \left( \textbf{prox}_{\zeta_k g_1}\left( \cdot\right)\right)$, and $B_{g_2} \in \partial_w^C\left(\textbf{prox}_{\zeta_k g_2}\left(\cdot\right)\right)$. Thus, using \eqref{eqn: Clarke subdifferential matrices} we can choose $B_{\delta,k_j},\ B_{g_1,k_j},\ B_{g_2,k_j}$ from the Bouligand subdifferential to improve computational efficiency, by reducing the active variables and constraint rows. To that end, we set $B_{\delta,k_j}$ as a diagonal matrix with
\begin{equation}\label{eqn: Clarke subdifferential of projection choice}
\left(B_{\delta,k_j}\right)_{(i,i)} \coloneqq \begin{cases} 
1, &\qquad \textnormal{if}\quad \left(\beta_k^{-1}z_{k} + x_{k_j}\right)_i \in \left(a_{l_i},a_{u_i}\right),\\
0, &\qquad  \textnormal{otherwise},
\end{cases},
\end{equation}
\noindent for all $i \in \{1,\ldots,n\}$, $B_{g_1,k_j}$ as the following diagonal matrix
\begin{equation}\label{eqn: Clarke subdifferential of prox choice}
\left(B_{g_1,k_j}\right)_{(i,i)} \coloneqq \begin{cases} 
1, &\qquad \textnormal{if}\quad \left| \left(\widehat{u}_{k_j}\right)_i\right| > \zeta_k D_{(i,i)},\ \textnormal{or}\quad D_{(i,i)} = 0,\\
0, &\qquad  \textnormal{otherwise},
\end{cases},
\end{equation}
\noindent for all $i \in \{1,\ldots,n\}$, where $\widehat{u}_{k_j} \coloneqq x_{k_j} - \zeta_k(r_{k}(x_{k_j},y_{k_j}))_{1:n}$, and $B_{g_2,k_j}$ as
\begin{equation}\label{eqn: Clarke subdifferential of prox g2 choice}
\left(B_{g_2,k_j}\right)_{(i,i)} \coloneqq \begin{cases} 
1, &\qquad \textnormal{if}\quad \left(w_{k_j}\right)_i - \zeta_k  \left(y_{k_j}\right)_i \leq 0,\ \textnormal{or}\ \left(w_{k_j}\right)_i - \zeta_k \left(y_{k_j}\right)_i \geq \zeta_k,\\
0, &\qquad  \textnormal{otherwise}
\end{cases},
\end{equation}
\noindent for all  $i \in \{1,\ldots,l\}$. 
\par Given the aforementioned choices for the projection matrices, we can now explicitly eliminate certain variables from the system in \eqref{Primal-dual SSN system}, in order to obtain a saddle-point system. To that end, from the second block-equation in \eqref{Primal-dual SSN system}, we have 
\[\left(I_l - B_{g_2,k_j}\right)d_w + \zeta_k B_{g_2,k_j} (d_y)_{1:l} = -\left(w_{k_j} - \textbf{prox}_{\zeta_k g_2} \left(w_{k_j} - \zeta_k (y_{k_j})_{1:l} \right) \right).\]
\noindent Let $\mathcal{B}_{g_2,k_j} \coloneqq \left\{i \in \{1,\ldots, l\} \colon \left(B_{g_2,k_j}\right)_{(i,i)} = 1\right\}$. Then, we obtain
\begin{equation} \label{eqn: dy parts getting constant}
(d_y)_{\mathcal{B}_{g_2,k_j}} = -\zeta_k^{-1}\left( w_{k_j} - \textbf{prox}_{\zeta_k g_2} \left(w_{k_j} - \zeta_k(y_{k_j})_{1:l} \right)\right)_{\mathcal{B}_{g_2,k_j}},
\end{equation}
\noindent where the right-hand side can be evaluated easily. Letting $\mathcal{N}_{g_2,k_j} \coloneqq \{1,\ldots,l\}\setminus \mathcal{B}_{g_2,k_j}$, we have
\begin{equation} \label{eqn: dw parts getting constant}
(d_w)_{\mathcal{N}_{g_2,k_j}} = - (w_{k_j})_{\mathcal{N}_{g_2,k_j}}.
\end{equation}
\noindent On the other hand, from the third block-equation of \eqref{Primal-dual SSN system}, we observe that
\[(d_w)_{\mathcal{B}_{g_2,k_j}} = -\left(w_{k_j} - d - C\left(x_{k_j}+d_x\right) - \beta_k^{-1}\left(y_{k_j} + d_y - y_k\right)_{1:l} \right)_{\mathcal{B}_{g_2,k_j}},\]
\noindent which can be computed after solving the reduced system. We define the following index sets
 \[\mathcal{B}_{g_1,k_j} \coloneqq \left \{ i \in \{1,\ldots,n\} \colon \left(B_{g_1,k_j}\right)_{(i,i)} = 1 \right\},\qquad \mathcal{N}_{g_1,k_j} \coloneqq \{1,\ldots,n\}\setminus \mathcal{B}_{g_1,k_j}.\]
\noindent From the first block equation of \eqref{Primal-dual SSN system}, we obtain
\[(d_x)_{\mathcal{N}_{g_1,k_j}} = -\left(x_{k_j} - \textbf{prox}_{\zeta_k g_1}\left(x_{k_j} - \zeta_k\left(r_{\beta_k,\rho_k}(x_{k_j},y_{k_j})\right)_{1:n} \right)\right)_{\mathcal{N}_{g_1,k_j}}. \]
\noindent After a straightforward (if tedious) calculation, we can pivot $(d_y)_{\mathcal{B}_{g_2,k_j}}$, $d_w$, and $(d_x)_{\mathcal{N}_{g_1,k_j}}$ in system \eqref{Primal-dual SSN system}. This results in the following saddle-point system:
\begin{equation} \label{eqn: detailed explicit SSN linear system}
\underbrace{\begin{bmatrix}- \left(\zeta_k^{-1}H_{k_j}\right)_{\left(\mathcal{B}_{g_1,k_j},\mathcal{B}_{g_1,k_j}\right)} & \begin{bmatrix} \widehat{C} \\ \widehat{A} 
\end{bmatrix}^\top\\
\begin{bmatrix} \widehat{C}  \\ \widehat{A} 
\end{bmatrix}  & \beta_k^{-1} I_{m+ \lvert\mathcal{N}_{g_2,k_j} \rvert}
\end{bmatrix}}_{\widehat{M}_{k_j}} \begin{bmatrix}
d_{x,\mathcal{B}_{g_1,k_j}}\\
\begin{bmatrix} \left(d_y\right)_{\mathcal{N}_{g_2,k_j}} \\ (d_y)_{(l+1:l+m)}
\end{bmatrix}
\end{bmatrix} = \widehat{r}_{k_j},
\end{equation}
\noindent where $\widehat{C} \coloneqq C_{\left({\mathcal{N}_{g_2,k_j}},\mathcal{B}_{g_1,k_j}\right)}, \ \widehat{A} \coloneqq A_{\mathcal{B}_{g_1,k_j}}$, and
\begin{equation*}
\begin{split}
\widehat{r}_{k_j} \coloneqq &\ \begin{bmatrix}
\zeta_k^{-1}\left(\widehat{M}_k(x_{k_j},w_{k_j},y_{k_j})\right)_{\mathcal{B}_{g_1,k_j}} \\
 \left(\left(-\widehat{M}_k(x_{k_j},w_{k_j},y_{k_j})\right)_{(n+l+1:n+2l)}\right)_{\mathcal{N}_{g_2,k_j}} 
\\
 \left(-\widehat{M}_k(x_{k_j},w_{k_j},y_{k_j})\right)_{(n+2l+1:n+2l+m)} 
\end{bmatrix}  \\
&\qquad + \begin{bmatrix}
  Q_{\left(\mathcal{B}_{g_1,k_j},\mathcal{N}_{g_1,k_j}\right)} d_{x,\mathcal{N}_{g_1,k_j}} + \left(C_{\left({\mathcal{B}_{g_2,k_j}},\mathcal{B}_{g_1,k_j}\right)} \right)^\top (d_y)_{\mathcal{B}_{g_2,k_j}}\\
  - C_{\left(\mathcal{N}_{g_2,k_j},\mathcal{N}_{g_1,k_j}\right)} (d_x)_{\mathcal{N}_{g_1,k_j}} + (d_w)_{\mathcal{N}_{g_2,k_j}} \\
  - A_{\mathcal{N}_{g_1,k_j}} d_{x,\mathcal{N}_{g_1,k_j}}
\end{bmatrix}.
\end{split}
\end{equation*}
\noindent Notice that the coefficient matrix of \eqref{eqn: detailed explicit SSN linear system} is symmetric and \emph{quasi-definite} (see \cite{SIAMOpt:Vander}). As such, it is strongly factorizable, that is, each symmetric permutation of it admits an $L\Delta L^\top$ decomposition, with $\Delta$ diagonal. Alternatively, one could eliminate further variables in order to obtain a positive definite linear system. This could be beneficial in certain applications, but is omitted here for the sake of generality. We note that $\widehat{C}$ contains only a subset of the columns and rows of $C$. Similarly, $\widehat{A}$ contains only a subset of the columns of $A$. As we will verify in practice, relatively close to the solution of \eqref{primal problem} the active-set matrices $\widehat{C}$, $\widehat{A}$, can be significantly smaller than $C$ and $A$, respectively, allowing the proposed approach to solve large-scale instances, without requiring excessive memory. Finally, we note that if we had included a proximal term for the auxiliary variables $w$ (i.e. if we were to directly apply the algorithm given in \cite{arxiv:PougkGondzio}), then $\widehat{C}$ would only contain a subset of the columns of $C$ but all of its rows (which, in general, are expected to be much more problematic, since it is often the case that $l \gg n$).

\subsection{Extension to arbitrary piecewise-linear terms} \label{subsection: extension to robust optimization}
\par Before closing this section we would like to notice that there are certain piecewise-linear terms that are not captured by the model in \eqref{primal problem}. In particular, given a set of $K$ pairs $(C_r,d_r) \in \mathbb{R}^{l\times n} \times \mathbb{R}^l$ (with $K\geq 3$), where $r \in \{1,\ldots,K\}$, we could consider problems of the following form
\begin{equation*} 
\begin{split}
\underset{x \in \mathbb{R}^n}{\text{min}}&\quad \left\{c^\top x + \frac{1}{2}x^\top Q x + \sum_{i = 1}^l \max_{r \in \{1,\ldots,R\}} \{\left(C_{r} x + d_{r}\right)_i\} + \|Dx\|_1 + \delta_{\mathcal{K}}(x)\right\}, \\ 
\textnormal{s.t.} &\quad \ A x = b.
\end{split}
\end{equation*}
\noindent In order to do so, we would have to introduce $K$ auxiliary vectors $w_r \in \mathbb{R}^l$, $r \in \{1,\ldots,K\}$, and reformulate the problem in the following form
\begin{equation*} 
\begin{split}
\underset{x \in \mathbb{R}^n}{\text{min}}&\quad \left\{c^\top x + \frac{1}{2}x^\top Q x + \sum_{i = 1}^l \max_{r \in \{1,\ldots,K\}} \{(w_r)_i\} + \|Dx\|_1 + \delta_{\mathcal{K}}(x)\right\}, \\ 
\textnormal{s.t.} &\quad \ C_{r} x + d_r - w_r = 0_l,\qquad  r \in \{1,\ldots,K\},\\
&\quad\ Ax = b.
\end{split}
\end{equation*}
\noindent Subsequently, in the derivation of the semismooth Newton method, we would have to evaluate the proximity operator of $\tilde{g}(u) \coloneqq \max_{i \in \{1,\ldots,K\}}(u_i)$, $u \in \mathbb{R}^K$, which admits a closed-form solution:
\[\left(\textbf{prox}_{\zeta \tilde{g}}(u)\right)_i = \min\{u_i,s\},\quad \textnormal{where } s \in \mathbb{R}\textnormal{ is such that }\sum_{i = 1}^K \left(\lvert u_i\rvert-s\right)_+ = \zeta.\]
\noindent Then, the Clarke derivative of the previous could be computed by utilizing \cite[Theorem 14.7]{arXiv:ClasValk}. This was not considered in this work for brevity of presentation, and is left as a future research direction. Indeed, this extension paves the way for the generalization of the proposed active-set scheme to a plethora of robust optimization problems that appear in practice (e.g. see \cite{OperRes:Berts}), as well as delivering an alternative to standard cutting-plane methods (e.g. see \cite{JSIAM:Kelley}), decomposition methods (e.g. see \cite{NumMath:Benders,OperRes:Birge,Econometrica:Dantzig,MathProg:Rusz}), or specialized interior point (e.g. see \cite{SIAMOpt:Gor}) approaches appearing in the literature.
\section{Warm-starting} \label{sec: warm start}
\par Following the developments in \cite{SIAMOpt:Lietal,arxiv:PougkGondzio}, we would like to find a starting point for Algorithm \ref{PMM algorithm} that is relatively close to a primal-dual solution. Indeed, this is crucial since, on the one hand, then we can expect to observe early linear convergence of Algorithm \ref{PMM algorithm}, while on the other hand, we can expect to be close to identifying the correct active-sets which in turn implies that the memory requirements of Algorithm \ref{primal-dual SNM algorithm} are significantly reduced. To that end, we employ a proximal alternating direction method of multipliers (pADMM; e.g. see \cite{SciComp:DengYin}). We reformulate \eqref{primal problem} by introducing an artificial variable $u \in \mathbb{R}^{n+l}$, as
\begin{equation} \label{primal problem ADMM reformulation} \tag{P'}
\begin{split}
\underset{(x,w,u)\ \in\ \mathbb{R}^n\times\mathbb{R}^l\times \mathbb{R}^{n+l}}{\text{min}} & \left\{ c^\top x  + \frac{1}{2}x^\top Q x +  \sum_{i = n+1}^{n+l} \left(u_{i}\right)_+ + \|D(u_{1:n})\|_1 + \delta_{\mathcal{K}}(u_{1:n})\right\},\\
\textnormal{s.t.}\qquad\ \ &\ \underbrace{\begin{bmatrix} C& -I_{l} & 0_{l,l+n}\\ A & 0_{m,l} & 0_{m,l+n} \\  \multicolumn{2}{c}{-I_{l+n}} & I_{l+n} \end{bmatrix}}_{M_{r}} \begin{bmatrix} x \\ w \\u\end{bmatrix} = \begin{bmatrix} -d \\ b\\ 0_{l+n}\end{bmatrix}.\\
\end{split}
\end{equation}
\par Given a penalty $\sigma > 0$, we associate the following augmented Lagrangian to \eqref{primal problem ADMM reformulation}
\begin{equation*}
\begin{split}
\widehat{\mathcal{L}}_{\sigma}(x,w,u,y) \coloneqq &\  c^\top  x  + \frac{1}{2}x^\top Q x +  \sum_{i = n+1}^{n+l} \left(u_{i}\right)_+ + \|D(u_{1:n})\|_1 + \delta_{\mathcal{K}}(u_{1:n})\\&\ - y^\top \left(M_r \begin{bmatrix} x \\ w \\u\end{bmatrix} - \begin{bmatrix} -d \\ b\\ 0_{l+n}\end{bmatrix} \right) + \frac{\sigma}{2}\left\|M_r \begin{bmatrix} x \\ w \\u\end{bmatrix}-\begin{bmatrix} -d \\ b\\ 0_{l+n}\end{bmatrix}\right\|^2,
\end{split}
\end{equation*}
\noindent where $y \in \mathbb{R}^{m+n+2l}$ denotes the dual multipliers associated to the linear equality constraints of  \eqref{primal problem ADMM reformulation}. Let an arbitrary positive definite matrix $R \in \mathbb{R}^{(n+l)\times(n+l)}$ be given, and assume the notation $\|v\|^2_{R} = v^\top R v$, for any $v \in \mathbb{R}^{n+l}$. Also, as in Section \ref{subsec: AL penalties}, we denote $g(x,w) = \|D x\|_1 + \sum_{i=1}^l \left( w_i\right)_+$.  Algorithm \ref{proximal ADMM algorithm} describes a proximal ADMM for the approximate solution of \eqref{primal problem ADMM reformulation}.
\renewcommand{\thealgorithm}{pADMM}
\begin{algorithm}[!ht]
\caption{(\emph{proximal ADMM})}
    \label{proximal ADMM algorithm}
    \textbf{Input:}  $\sigma > 0$, $R \succ 0$, $\gamma \in \left(0,\frac{1+\sqrt{5}}{2}\right)$, $(x_0,w_0,u_0,y_{0}) \in \mathbb{R}^{2n + 3l + m}$.
\begin{algorithmic}
\For {($k = 0,1,2,\ldots$)}
\begin{flalign*}
u_{k+1} &= \underset{u}{\arg\min}\left\{\widehat{\mathcal{L}}_{\sigma}\left(x_k,w_k,u,y_{k}\right) \right\}\equiv \Pi_{\mathcal{K}\times \mathbb{R}^l}\left(\textbf{prox}_{\sigma^{-1} g}\left(\begin{bmatrix} x_k\\ w_k\end{bmatrix} + \sigma^{-1} (y_{k})_{(m+l+1:m+n+2l)}\right)\right).\\ 
\begin{bmatrix} x_{k+1} \\ w_{k+1} \end{bmatrix} &= \underset{x,\ w}{\arg\min}\left\{\widehat{\mathcal{L}}_{\sigma}\left(x,w,u_{k+1},y_{k}\right) + \frac{1}{2}\left\|\begin{bmatrix} x-x_k \\ w - w_k\end{bmatrix}\right\|_{R}^2\right\}.\\ 
y_{k+1} &= y_{k} - \gamma\sigma\left(M_r \begin{bmatrix} x \\ w \\u\end{bmatrix} - \begin{bmatrix} -d \\ b\\ 0_{l+n}\end{bmatrix} \right).
\end{flalign*}
\EndFor
\end{algorithmic}
\end{algorithm}
\par Let us notice that under certain standard assumptions on \eqref{primal problem ADMM reformulation}, Algorithm \ref{proximal ADMM algorithm} converges globally, potentially at a linear rate (see \cite{SciComp:DengYin}). We can employ the regularization matrix $R$ as a means of ensuring memory efficiency of Algorithm \ref{proximal ADMM algorithm}. For example, we can recover the \emph{prox-linear ADMM} \cite[Section 1.1]{SciComp:DengYin}, where given some sufficiently large constant $\hat{\sigma} > 0$, one defines
\[ R \coloneqq \hat{\sigma}I_{n+l}  - \sigma\begin{bmatrix} C^\top C + A^\top  A +\sigma^{-1}\textnormal{Off}(Q)& -C^\top \\ - C & 0_{l,l} \end{bmatrix}\succ 0.\]
\noindent  Given this choice of $R$, the second step of Algorithm \ref{proximal ADMM algorithm} consists of only matrix-vector multiplications with $A,\ C$ and $Q$, and thus no (non-diagonal) matrix inversion is required. If memory is not an issue, one could use a positive-definite diagonal regularizer $R$, yielding a standard regularized ADMM. In our implementation, the user can specify which of the two strategies should be employed. In either case, the first step of Algorithm \ref{proximal ADMM algorithm} is trivial to solve, since we know that the proximity operator of $g(\cdot)$ can be computed analytically as in Section \ref{subsec: SSN method}. 
\par Finally, once an approximate solution $(\tilde{x},\tilde{w},\tilde{u},\tilde{y})$ is computed, we set the starting point of Algorithm \ref{PMM algorithm} as $(x_0,w_0,y_0,z_0) = \left(\tilde{x},\tilde{w},\left(\tilde{y}\right)_{(1:m+l)},\tilde{z}\right)$, where
\[\tilde{z} = \left(\tilde{y}\right)_{(m+l+1:m+l+n)} - \Pi_{\partial g_1(\tilde{u}_{1:n})}\left(\left(\tilde{y}\right)_{(m+l+1:m+l+n)}\right).\]
\noindent Indeed, an optimal primal-dual solution of \eqref{primal problem ADMM reformulation} satisfies 
\[\left(\tilde{y}^*\right)_{(m+l+1:m+2l+n)} \in \partial g\left(\tilde{u}^*\right) + \partial \delta_{\mathcal{K}\times \mathbb{R}^l}\left(\tilde{u}^*\right),\]
\noindent thus the characterization of $\tilde{z}$ in Algorithm \ref{PMM algorithm} follows from the Appendix, where one can also find the termination criteria of Algorithm \ref{proximal ADMM algorithm} (see \eqref{termination criteria for pADMM}).
\section{Applications and numerical results} \label{sec: Applications}
\par In this section we present various applications that can be modeled by problem \eqref{primal problem}, focusing on portfolio optimization, quantile regression, and binary classification. In particular, we first discuss (and numerically demonstrate) the effectiveness of the approach for the solution of single-period mean-risk portfolio optimization problems, where risk is measured via the \emph{conditional value at risk} or the \emph{mean absolute semi-deviation}. Subsequently, we apply the proposed scheme for the solution of quantile regression problems, demonstrating its scalability. Finally, we apply the active-set scheme for the solution of binary classification problems via linear support vector machines on certain large-scale datasets. 
\paragraph{Implementation details} Before proceeding to the numerical results, we mention certain implementation details of the proposed algorithm. The implementation follows closely the developments in \cite{arxiv:PougkGondzio}, and can be found on GitHub\footnote{\url{https://github.com/spougkakiotis/Active_set_method_for_CQP_piecewise_LP}}. The code is written in MATLAB, and the experiments are run on a PC with a 2.2GHz Intel core i7-8750H processor (hexa-core), 16GB of RAM, using the Windows 10 operating system.
\par We run Algorithm \ref{proximal ADMM algorithm} (warm-start) for at most 100 iterations (or until a 3-digit accurate solution is found). The user can choose whether the warm-starting scheme should be matrix-free or not. In the presented experiments, the matrix-free scheme was only used for the largest quantile regression and classification problems, for which standard ADMM crashed due to excessive memory requirements. Then, starting with $\beta_0 = 10$, $\rho_0 = 50$, we run Algorithms \ref{PMM algorithm}-\ref{primal-dual SNM algorithm}. Following \cite{arxiv:PougkGondzio}, when solving the PMM sub-problems using Algorithm \ref{primal-dual SNM algorithm} we use a predictor-corrector-like heuristic in which the first iteration is accepted without line-search and then line-search is activated for subsequent iterations. Algorithm \ref{PMM algorithm} is allowed to run for at most 200 iterations, while Algorithm \ref{primal-dual SNM algorithm} is stopped after 20 inner iterations. An iterate is accepted as optimal if the conditions given in \eqref{termination criteria for PD-PMM} are satisfied for the tolerance specified by the user. Most other implementation details follow directly from the developments in Sections \ref{subsec: AL penalties}--\ref{subsec: SSN method}. We refer the reader to the implementation on GitHub for additional details.
\par In the presented experiments we compare the proposed approach against the interior point solver given in \cite{COAP:PougkGond}, the implementation of which can be found on GitHub\footnote{\url{https://github.com/spougkakiotis/IP_PMM}}, as well as the ADMM-based OSQP solver given in \cite{osqp}, the implementation of which can also be found on GitHub\footnote{\url{https://github.com/osqp/osqp-matlab}}. We note that most problems considered herein were especially challenging for OSQP (mostly due to requesting highly accurate solutions), and for that reason we allow it to run for 50,000 iterations (unlike its default iteration threshold, which is 4000).

\subsection{Portfolio optimization}
\par We first consider the mean-risk portfolio selection problem (originally proposed in \cite{JoFinance:Markowitz}), where we minimize some convex risk measure, while keeping the expected return of the portfolio above some desirable level. A variety of models for this problem have been intensely analyzed and solved in the literature (e.g. see \cite{JBanFin:AlexColYi,JRisk:KroUryPal,JRisk:RockUry}). The departure from the variance as a measure of risk often allows for great flexibility in the decision making of investors, enabling them to follow potential regulations as well as to better control the risk associated with an ``optimal" portfolio. Optimality conditions and existence of solutions for several general deviation measures have been characterized in \cite{MathProg:RockUryZaba}. The comparison of portfolios obtained by minimizing different risk measures has also been considered (e.g. see \cite{ssnr:Rametal,FinResLet:RigBor,Prod:Silva_etal}). 
\par The method presented in this paper is general and not related to a specific model choice. Indeed, we would like to showcase the efficiency of our approach for obtaining accurate solutions to portfolio optimization problems with various risk measures of practical interest. Thus, we focus on the solution of a standard portfolio selection model that has been used in the literature. All numerical results are obtained on real-world datasets. As a result the problems are of medium-scale. Nonetheless, even for such medium-scale problems, we will be able to demonstrate the efficiency of the proposed approach, when compared to the efficient interior point method employed in the literature \cite{SIREV:DeSimone_etal} for similar problems. We also note that both second-order solvers (i.e. IPM and active-set) significantly outperform OSQP on these instances, however it was included in the comparison for completeness. Some large-scale instances will be tackled in the following subsections, where the method will be applied to quantile regression and binary classification instances.
\par Let $x \in \mathbb{R}^n$ represent a portfolio of $n$ financial instruments, such that 
\[x_i \in \left[a_{l_i},a_{u_i}\right],\quad  \textnormal{with}\ a_{l_i} \geq -1,\ a_{u_i} \leq 1,\ \textnormal{for all}\ i = 1,\ldots,n,\quad \textnormal{and}\ \sum_{i=1}^n x_i = 1. \]
\noindent This requirement indicates that short positions for each stock are restricted by some percentage ($a_{l_i}\%$) of the available wealth (assuming $a_{l_i} < 0$), and no more than $a_{u_i}\%$ of the total wealth can be invested to instrument $i$. Let $\bm{\xi} \in \mathbb{R}^n$ denote a random vector, the $i$-th entry of which represents the random return of the $i$-th instrument. Then, the random loss (i.e. the negative of the random return) of a given portfolio $x$ is given by $f(x,\bm{\xi}) \coloneqq -x^\top \bm{\xi}$. In this paper we assume that $\bm{\xi}$ follows some continuous distribution $p(\bm{\xi})$, as well as that there is a one-to-one correspondence between percentage return and monetary value (as in \cite[Section 3]{JRisk:RockUry}).  Additionally, given some expected benchmark return $r$ (e.g. the \emph{market index}), we only consider portfolios that yield an expected return above a certain threshold, i.e. $\mathbb{E} [-f(x,\bm{\xi})] \geq r$.
\par Finally, given the previously stated constraints, we would like to minimize some convex risk measure $\varrho(\cdot)$ of interest. However, in order to make sure that the problem is well-posed, while the transaction costs are not excessive, we include an $\ell_1$ term in the objective. This is a well-known modeling choice that yields sparse portfolios, and thus regulates the transaction costs in the single-period portfolio setting (e.g. see \cite{JBanFin:AlexColYi,AnnOpRes:Corsaro_etal,JRisk:KroUryPal}). Additionally, with an appropriate tuning of the $\ell_1$ regularization parameter $\tau > 0$, one could also control the amount of short positions (see \cite{AnnOpRes:Corsaro_etal}). It should be mentioned here that in the multi-period setting such an $\ell_1$ term does not guarantee a reduction in the transaction costs (but mostly in the \emph{holding costs}), and an additional \emph{total variation} term should also be added in the objective (see \cite{AnnOpRes:Corsaro_etat2,SIREV:DeSimone_etal}). This is omitted here, since we focus on the single-period case, but the model in \eqref{primal problem} could easily incorporate such an additional term (see Remark \ref{remark: l1 norm reformulation}). By putting everything together, the model reads as
\begin{equation} \label{model: vanilla portfolio selection}
\begin{split}
\min_{x} &\ \varrho\left(f(x,\bm{\xi})\right) + \tau \|x\|_1,\\
\textnormal{s.t.}&\ \sum_{i=1}^n x_i = 1,\\
&\ \mathbb{E}\left[-f(x,\bm{\xi})\right] \geq r,\\
&\ x_i \in \left[a_{l_i},a_{u_i}\right],\qquad i = 1,\ldots,n.
\end{split}
\end{equation}
\par There are several methods for solving such stochastic problems; let us mention two important variants. There is the \emph{parametric approach} (e.g. as in \cite{JBanFin:AlexColYi,JRisk:RockUry}), where one assumes that the returns follow some known distribution which is subsequently sampled to yield finite-dimensional optimization problems, and the \emph{sampling approach} (e.g. as in \cite{JRisk:KroUryPal}), where one obtains a finite number of samples (without assuming a specific distribution). Such samples are often obtained by historical observations, and this approach is also followed in this paper. It is well-known that historical data cannot fully predict the future (see \cite{JRisk:KroUryPal}), however it is a widely-used practice. The reader is referred to \cite{EFM:Jorion} for an extensive discussion onprobabilistic models for portfolio selection problems.
\par Additional soft or hard constraints can be included when solving a portfolio selection problem. Such constraints can either be incorporated directly via the use of a model (e.g. see \cite[Section 2]{NowPublishers:Boyd_etal}) as hard constraints or by including appropriate penalty terms in the objective (soft constraints). It is important to note that the model given in \eqref{primal problem} is quite general and as a result has great expressive power, allowing one to incorporate various important risk measures (and their combinations), as well as modeling constraints of interest.
\paragraph{Real-world datasets:} In what follows, we solve two different instances of problem  \eqref{model: vanilla portfolio selection}. In particular, we consider two potential risk measures; the conditional value at risk (e.g. see \cite{JRisk:RockUry}), as well as the mean absolute semi-deviation (e.g. see \cite[Section 6.2.2.]{SIAM:Shapiro}, noting that this is in fact equivalent to the mean absolute deviation originally proposed in \cite{ManagSci:KonnoYama}). In each of the above cases problem \eqref{model: vanilla portfolio selection} has the form of \eqref{primal problem} and thus Algorithm \ref{PMM algorithm} can directly be applied. 
\par We showcase the effectiveness of the proposed approach on 6 real datasets taken from \cite{DataBrief:Bruni_etal}. Each dataset contains time series for weekly asset returns and market indexes for  different major stock markets, namely, DowJones, NASDAQ100, FTSE100, SP500, NASDAQComp, and FF49Industries. In the first 5 markets the authors in \cite{DataBrief:Bruni_etal} provide the market indexes, while for the last dataset the uniform allocation strategy is considered as a benchmark. Additional information on the datasets is collected in Table \ref{Table: portfolio optimization datasets}. We note that stocks with less than 10 years of observations have been disregarded.
\begin{table}[!ht]
\centering
\caption{Portfolio optimization datasets.\label{Table: portfolio optimization datasets}}
\scalebox{1}{
\begin{tabular}{llll}     
\specialrule{.2em}{.05em}{.05em} 
    \textbf{Name} & \textbf{\# of assets} & \textbf{\# of data points} & \textbf{Timeline} \\ \specialrule{.1em}{.3em}{.3em} 
DowJones & 28 & 1363 & Feb. 1990--Apr. 2016\\ 
NASDAQ100 & 82 & 596 & Nov. 2004--Apr. 2016\\ 
FTSE100 & 83 & 717 & Jul. 2002--Apr. 2016\\ 
FF49Industries & 49 & 2325 & Jul. 1969--Jul. 2015\\ 
SP500 & 442 & 595 & Nov. 2004--Apr. 2016\\ 
NASDAQComp & 1203 & 685 & Feb. 2003--Apr. 2016\\ 
\specialrule{.2em}{.05em}{.05em} 
\end{tabular}}
\end{table}

\subsubsection{Conditional value at risk}
\par First, we consider portfolio optimization problems that seek a solution minimizing the conditional value at risk; a measure which is known to be \emph{coherent} (see \cite{MathFinance:Artzner} for a definition of coherent risk measures). In particular using the notation introduced earlier, we consider the following optimization problem
\begin{equation} \label{CVaR original problem}
\underset{x \in \mathbb{R}^n}{\text{min}} \left\{  \textnormal{CVaR}_{\alpha}\left(f(x,\bm{\xi})\right) + \tau\|x\|_1 + \delta_{\mathcal{K}}(x)\right\}, \qquad \textnormal{s.t.}\ A x = b,
\end{equation}
\noindent where $f(x,\bm{\xi})$ is the random cost function, $A \in \mathbb{R}^{m\times n}$ models the linear constraint matrix of problem \eqref{model: vanilla portfolio selection} (where an auxiliary variable has been introduced to transform the inequality constraint involving the expectation into an equality), and $\mathcal{K} \coloneqq [a_l,a_u]$. In the above $1-\alpha \in (0,1)$ is the confidence level. It is well-known (\cite{JRisk:RockUry}) that given a continuous random variable $X$, the conditional value at risk can be computed as
\begin{equation*} 
\textnormal{CVaR}_{\alpha}\left(X\right)\coloneqq \min_{t \in \mathbb{R}}\left\{t + \alpha^{-1}\mathbb{E}\left[ \left(X - t\right)_+ \right]\right\}.
\end{equation*}
\noindent  We can write problem \eqref{CVaR original problem} in the following equivalent form:
\begin{equation*}
\underset{(x,t) \in \mathbb{R}^n\times \mathbb{R}}{\text{min}} \left\{ t + \frac{1}{l \alpha} \sum_{i = 1}^l \left(-\xi_i^\top x - t\right)_+ + \|Dx\|_1 + \delta_{\mathcal{K}}(x)\right\},\qquad \textnormal{s.t.}\ Ax = b,
\end{equation*}
\noindent where the expectation has been substituted by summation  since we assume the availability of a dataset $\{\xi_1,\ldots,\xi_l\}$. Introducing an auxiliary variable $w \in \mathbb{R}^l$, the previous can be re-written as
\begin{equation*} 
\begin{split}
\underset{(x,t,w) \in \mathbb{R}^n\times \mathbb{R}\times\mathbb{R}^l}{\text{min}} & \left\{ t +  \sum_{i = 1}^l \left(w_i\right)_+ + \|Dx\|_1 + \delta_{\mathcal{K}}(x)\right\},\\
\textnormal{s.t.}\qquad\ \ &\ \frac{1}{l \alpha}\left(-\xi_i^\top x-t\right) - w_i = 0, \quad i = 1,\ldots,l,\\
&\ Ax = b.
\end{split}
\end{equation*}
\par We solve problem \eqref{CVaR original problem} using the proposed active-set method (AS), the \emph{interior point-proximal method of multipliers} (IP-PMM) given in \cite{COAP:PougkGond}, and OSQP (\cite{osqp}). We allow any short position, that is $a_{l_i} = -1$, and restrict investing more that $60\%$ of the available wealth on a single stock (i.e. $a_{u_i} = 0.6$). We set $\texttt{tol} = 10^{-5}$ and $\tau = 10^{-2}$, and run the three methods for each of the datasets described in Table \ref{Table: portfolio optimization datasets} for varying confidence level. We report the confidence parameter $\alpha$, the number of PMM, SSN, IP-PMM and OSQP iterations, the CPU time needed by each of the three schemes, as well as the number of factorizations used within the active set (PMM-SSN) scheme. Indeed, it often happens that the active-set is not altered from one iteration to the next, and the factorization does not need to be re-computed. The results are collected in Table \ref{Table CVaR portfolio selection: varying confidence level}. In all the numerical results that follow, the lowest running time exhibited by a solver, assuming it successfully converged, is presented in bold. 
\begin{table}[!ht]
\caption{CVaR portfolio selection: varying confidence level (\texttt{tol} = $10^{-5}$, $\tau = 10^{-2}$).\label{Table CVaR portfolio selection: varying confidence level}}
\centering
\scalebox{0.95}{\begin{threeparttable}
\begin{tabular}{llllllll}     
\specialrule{.2em}{.05em}{.05em} 
    \multirow{2}{*}{\textbf{Dataset}} & \multirow{2}{*}{$\bm{\alpha}$} & \multicolumn{3}{c}{\textbf{Iterations}}     & \multicolumn{3}{c}{\textbf{Time (s)}}  \\   \cmidrule(l{2pt}r{2pt}){3-5} \cmidrule(l{2pt}r{2pt}){6-8}
&   &  {{PMM}(SSN)[Fact.]} & {{IP--PMM}} & {{OSQP}} &  {{AS}} & {{IP--PMM}}  & {{OSQP}}   \\ \specialrule{.1em}{.3em}{.3em} 
 \multirow{3}{*}{{DowJones}} & $0.05$ & 35(169)[142] & 18 &   32,575& \textbf{0.65}  & 0.71 & 10.94   \\
    & $0.10$ & 36(166)[145] & 21 & 42,625 & \textbf{0.75} &  0.83 & 13.61 \\
        & $0.15$ & 36(144)[125] & 13 & 43,100 & 0.67  & \textbf{0.53} & 14.45 \\ \specialrule{.00002em}{.1em}{.1em}
 \multirow{3}{*}{NASDAQ100} & $0.05$ & 31(156)[138] & 19 & 30,575 & 0.66 & \textbf{0.60} & 7.92   \\
     & $0.10$ & 33(168)[150] & 16 & 34,700 & 0.65 & \textbf{0.50} & 8.81 \\
    &  $0.15$ & 33(163)[142] & 15 & 40,500 & 0.62 & \textbf{0.57}  & 10.46\\ \specialrule{.00002em}{.1em}{.1em}
\multirow{3}{*}{FTSE100} & $0.05$ & 32(179)[164] & 17 &   21,825& \textbf{0.77}  & 0.90 & 10.96   \\
    & $0.10$ & 34(171)[147] & 19 & 35,825 & \textbf{0.80} &  0.99 & 17.70 \\
        & $0.15$ & 35(179)[156] & 11 & 36,550 & 0.93  & \textbf{0.59} & 18.24 \\ \specialrule{.00002em}{.1em}{.1em}
\multirow{3}{*}{FF49Industries} & $0.05$  & 41(215)[186] & 27 &  $50,000^{\ddagger}$\tnote{1} & \textbf{1.77}  & 4.65 & $59.16^{\ddagger}$  \\
     & $0.10$ & 42(224)[183] & 22 &  $50,000^{\ddagger}$ & \textbf{2.37}  & 3.62 & $57.92^{\ddagger}$\\
      & $0.15$ & 39(205)[157] & 11 & $50,000^{\ddagger}$ & 2.66  & \textbf{2.18} &  $51.25^{\ddagger}$\\   \specialrule{.00002em}{.1em}{.1em}
 \multirow{3}{*}{SP500} & $0.05$  & 30(203)[193] & 20  & 21,075 & \textbf{5.93}  & 12.82 & 65.70  \\
    & $0.10$ & 34(199)[187] & 11 & 32,375 & \textbf{5.88}  & 6.98 & 99.39 \\
      & $0.15$ & 34(173)[160] & 17 &  39,225& \textbf{5.78}  & 10.17 & 119.95\\  \specialrule{.00002em}{.1em}{.1em}
     \multirow{3}{*}{NASDAQComp} & $0.05$ & 28(198)[192] & 10 & 16,100 & \textbf{22.07}  & 35.52 & 143.82  \\
     & $0.10$ & 31(171)[169] & 14 & 31,275 & \textbf{23.31}  & 48.99 & 291.97 \\
      & $0.15$ & 34(173)[167] & 19 &  44,775& \textbf{21.67}  & 67.14 & 403.75\\
\specialrule{.2em}{.05em}{.05em} 
\end{tabular}\begin{tablenotes}
\item[1] $\ddagger$ indicates that the solver reached the maximum number of iterations.
\end{tablenotes}
\end{threeparttable}}
\end{table}
\par From Table \ref{Table CVaR portfolio selection: varying confidence level} we observe that for the smaller instances (DowJones, NaSDAQ100, FTSE100, FF49Industries) both second-order methods perform quite well and are comparable in terms of CPU time. Nonetheless, even in this case, the proposed active-set solver requires significantly less memory. Indeed, this can be seen by the fact that the two methods achieve similar times but the active-set scheme is performing significantly more factorizations. On the other hand, for the larger instances (SP500, NASDAQComp) the proposed active-set scheme outperforms the interior point method significantly. This is mostly due to the efficiency gained from the lower memory requirements and cheaper factorizations of the active-set solver. Also, we observe that both methods are quite robust with respect to the confidence level and consistently outperform OSQP, which struggles to find a 5-digit accurate solution. We note that OSQP could potentially be competitive for smaller tolerances (e.g. for finding a 3-digit accurate solution), however, the application under consideration dictates that an accurate solution is needed, since a small improvement in the portfolio output can translate into huge profits in practice.
\par In order to observe the behaviour of the three solvers under a different $\ell_1$ regularization value, we fix $\texttt{tol} = 10^{-5}$, $\tau = 10^{-1}$, and run them for varying confidence level. The results are collected in Table \ref{Table CVaR portfolio selection: varying confidence level 2}. Similar observations can be drawn from Table \ref{Table CVaR portfolio selection: varying confidence level 2}, while noting that both second-order approaches remain efficient for different values of the regularization parameter $\tau$. Notice, however, that one should be careful when choosing this regularization parameter, since otherwise the obtained portfolios could be meaningless. In light of this, we have chosen values for $\tau$ that yield reasonable portfolios, with controlled short positions. Additionally we should note that both second-order solvers were able to solve all problems in higher accuracy (e.g. up to 6- or 7-digits of accuracy), however, the differences in the obtained portfolios (in terms of positions and/or associated risk) were negligible and hence such accuracies were not considered here. 
\begin{table}[!ht]
\centering
\caption{CVaR portfolio selection: varying confidence level (\texttt{tol} = $10^{-5}$, $\tau = 10^{-1}$).\label{Table CVaR portfolio selection: varying confidence level 2}}
\scalebox{0.95}{
\begin{tabular}{llllllll}     
\specialrule{.2em}{.05em}{.05em} 
    \multirow{2}{*}{\textbf{Dataset}} & \multirow{2}{*}{$\bm{\alpha}$} & \multicolumn{3}{c}{\textbf{Iterations}}     & \multicolumn{3}{c}{\textbf{Time (s)}}  \\   \cmidrule(l{2pt}r{2pt}){3-5} \cmidrule(l{2pt}r{2pt}){6-8}
&   &  {{PMM}(SSN)[Fact.]} & {{IP--PMM}}  & {{OSQP}} &  {{AS}} & {{IP--PMM}}  &{{OSQP}}   \\ \specialrule{.1em}{.3em}{.3em} 
 \multirow{3}{*}{{DowJones}} & $0.05$ & 35(168)[144] & 21 &24,975  & \textbf{0.57} & 0.85& 8.31   \\
    & $0.10$ & 36(154)[133] & 20 & 21,925  & \textbf{0.64} &  0.79 & 7.32 \\
        & $0.15$ & 36(143)[128] & 17 & 29,625 & \textbf{0.68}  & 0.69 & 9.95 \\ \specialrule{.00002em}{.1em}{.1em}
 \multirow{3}{*}{NASDAQ100} & $0.05$ & 34(157)[124] & 19 & 17,225 & \textbf{0.46} & 0.57 & 4.50   \\
     & $0.10$ & 34(155)[133] & 16 & 18,100 & 0.59 & \textbf{0.49} & 4.80\\
    &  $0.15$ & 35(157)[132] & 16 & 20,775 & 0.56&  \textbf{0.48} & 5.49 \\ \specialrule{.00002em}{.1em}{.1em}
\multirow{3}{*}{FTSE100} & $0.05$  & 34(166)[137] & 21 & 18,900 & \textbf{0.65} & 1.17 & 10.10 \\
     & $0.10$ & 34(174)[157] & 16 & 24,750 &  0.92 & \textbf{0.86}  & 13.15 \\
      & $0.15$ & 36(166)[147] & 17 &  28,850 & 0.95  & \textbf{0.92} & 15.53\\ \specialrule{.00002em}{.1em}{.1em}
\multirow{3}{*}{FF49Industries} & $0.05$  & 43(215)[171] & 19 & 48,000 & \textbf{1.63}  & 3.00 & 59.33  \\
     & $0.10$ & 41(184)[146] & 22 & 40,475 & \textbf{1.57}  & 3.53 & 48.48\\
      & $0.15$ & 41(160)[136] & 19 & 35,525 & \textbf{1.45} & 3.06 & 40.91\\   \specialrule{.00002em}{.1em}{.1em}
 \multirow{3}{*}{SP500} & $0.05$  & 32(163)[144] & 20 & 21,425 & \textbf{4.68}  & 12.26 & 66.29  \\
    & $0.10$ & 33(165)[143] & 18 & 22,100 & \textbf{5.86}  & 12.43 & 73.56 \\
      & $0.15$ & 35(169)[148] & 19 & 32,225 & \textbf{5.61}  & 11.51 & 98.25\\  \specialrule{.00002em}{.1em}{.1em}
     \multirow{3}{*}{NASDAQComp} & $0.05$ & 32(190)[182] & 26 & 17,200 & \textbf{15.55}  & 90.78  & 157.83 \\
     & $0.10$ & 34(180)[172] & 15 &  13,925 & \textbf{22.75}  & 51.81 & 127.26 \\
      & $0.15$ & 34(170)[167] & 18 &  14,675& \textbf{24.69}  & 62.40 & 137.42\\
\specialrule{.2em}{.05em}{.05em} 
\end{tabular}}
\end{table}
\subsubsection{Mean absolute semi-deviation}
\par Next, we consider portfolio optimization problems that seek a solution minimizing the mean absolute semi-deviation, which is also coherent. We consider the following optimization problem
\begin{equation} \label{MAD original problem}
\underset{x \in \mathbb{R}^n}{\text{min}} \left\{  \textnormal{MAsD}\left(f(x,\bm{\xi})\right) + \tau\|x\|_1 + \delta_{\mathcal{K}}(x)\right\}, \qquad \textnormal{s.t.}\ A x = b,
\end{equation}
\noindent where, given a continuous random variable $X$, the associated risk is defined as 
\begin{equation*}
\textnormal{MAsD}\left(X\right)\coloneqq \mathbb{E}\left[ \left(X - \mathbb{E}[X]\right)_+\right] \equiv \frac{1}{2} \mathbb{E}\left|X- \mathbb{E}[X]\right|,
\end{equation*}
\noindent where the equivalence follows from \cite[Prospotion 6.1]{SIAM:Shapiro}. Given a dataset $\{\xi_1,\ldots,\xi_l\}$, problem \eqref{MAD original problem} can be written as
\begin{equation*} 
\begin{split}
\underset{(x,w) \in \mathbb{R}^n\times\mathbb{R}^l}{\text{min}} & \left\{ \sum_{i = 1}^l \left(w_i\right)_+ + \|Dx\|_1 + \delta_{\mathcal{K}}(x)\right\},\\
\textnormal{s.t.}\quad\ \ &\ \frac{1}{l}\left(-\xi_i^\top x + \bm{\mu}^\top x\right)- w_i = 0, \quad i = 1,\ldots,l,\\
&\ Ax = b,
\end{split}
\end{equation*}
\noindent where $\bm{\mu} \coloneqq \frac{1}{l}\sum_{i= 1}^l \xi_i^\top x$. Note that this model is in the form of \eqref{primal problem}.
\par We fix $\texttt{tol} = 10^{-5}$ and run the two methods on the 6 datasets for two different sensible values of the regularization parameter $\tau$. The results are collected in Table \ref{Table MAsD portfolio selection: varying confidence level}.
\begin{table}[!ht]
\caption{MAsD portfolio selection: varying regularization (\texttt{tol} = $10^{-5}$).\label{Table MAsD portfolio selection: varying confidence level}}
\centering
\scalebox{0.95}{\begin{threeparttable}
\begin{tabular}{llllllll}     
\specialrule{.2em}{.05em}{.05em} 
    \multirow{2}{*}{\textbf{Dataset}} & \multirow{2}{*}{$\bm{\tau}$} & \multicolumn{3}{c}{\textbf{Iterations}}     & \multicolumn{3}{c}{\textbf{Time (s)}}  \\   \cmidrule(l{2pt}r{2pt}){3-5} \cmidrule(l{2pt}r{2pt}){6-8}
&   &  {{PMM}(SSN)[Fact.]} & {{IP--PMM}} & {{OSQP}}  &  {{AS}} & {{IP--PMM}}  & {{OSQP}}   \\ \specialrule{.1em}{.3em}{.3em} 
 \multirow{2}{*}{{DowJones}} & $0.01$ & 51(142)[139] & 10 & 32,350 & 1.05 &  \textbf{0.56} & 10.47\\
        & $0.05$ & 52(142)[139] & 13  & 14,050 & 0.72 &  \textbf{0.54}  & 4.60\\ \specialrule{.00002em}{.1em}{.1em}
 \multirow{2}{*}{NASDAQ100} & $0.01$ & 47(167)[160] & 11 & 18,800 & 0.67 & \textbf{0.26} & 5.02 \\
    &  $0.05$ & 48(176)[161] & 10 & 40,550 &  0.67 & \textbf{0.27} & 10.72\\ \specialrule{.00002em}{.1em}{.1em}
\multirow{2}{*}{FTSE100} & $0.01$ & 46(153)[150] & 10 & 18,325 & 0.79 & \textbf{0.27} & 9.56 \\
      & $0.05$ & 46(153)[150] & 14 & 47,125 & 0.74  & \textbf{0.35} & 23.29\\ \specialrule{.00002em}{.1em}{.1em}
\multirow{2}{*}{FF49Industries} & $0.01$  & 53(141)[136] & 11 &  $50,000^{\ddagger}$\tnote{1} & \textbf{1.49}  & 1.72 & $57.59^{\ddagger}$  \\
      & $0.05$ & 52(137)[132] & 9 & 12,050& \textbf{1.35}  & 1.51 & 14.24\\   \specialrule{.00002em}{.1em}{.1em}
 \multirow{2}{*}{SP500} & $0.01$  & 42(165)[163] & 17 & 28,375 & 5.33  & \textbf{2.91} & 87.10  \\
      & $0.05$ & 41(157)[153] & 17 & 41,275 & 4.39  & \textbf{3.07} & 128.09\\  \specialrule{.00002em}{.1em}{.1em}
     \multirow{2}{*}{NASDAQComp} & $0.01$ & 44(143)[133] & 17 & 34,175 & \textbf{9.75}  & 9.91 & 310.12   \\
      & $0.05$ & 44(143)[132] & 18 & 41,000 &  10.22  & \textbf{9.72} & 367.32 \\
\specialrule{.2em}{.05em}{.05em} 
\end{tabular}
\begin{tablenotes}
\item[1] $\ddagger$ indicates that the solver reached the maximum number of iterations.
\end{tablenotes}
\end{threeparttable}}
\end{table}
\par From Table \ref{Table MAsD portfolio selection: varying confidence level} we observe that both second-order schemes are robust and comparable for all instances. In this case, the larger instances (SP500, NASDAQComp) were solved in comparable time by both solvers. Nevertheless, it is important to note that the active-set scheme is the better choice, since it has significantly less memory requirements. Again, this can be readily observed from the fact that its factorizations are significantly cheaper compared to those of IP-PMM (since the active-set scheme performs significantly more factorizations and still achieves comparable performance). Indeed, as will become clear when solving large-scale quantile regression and binary classification instances, the proposed active-set solver scales better than IP-PMM or OSQP. Additionally, we also observe that the method is robust (i.e. converges reliably to an accurate solution). Finally, as in the case of the CVaR instances, OSQP struggled to find accurate solutions at a reasonable number of iterations, making it a less efficient choice for such problems.

\subsubsection{Extensions and alternative risk measures}
\par Let us notice that the presented methodology can be easily extended to the multi-period case \cite{NowPublishers:Boyd_etal,MathFin:LiNg}. Then, one could include an additional fused-lasso term in the objective function in order to ensure low transaction costs. It is important to note that in this case the $\ell_1$ term added in the objective has the effect of reducing holding costs as well as short positions (e.g. see \cite{AnnOpRes:Corsaro_etat2,AnnOpRes:Corsaro_etal,SIREV:DeSimone_etal}). As noted in Remark \ref{remark: l1 norm reformulation}, the additional fused-lasso term can be easily incorporated in the objective of \eqref{primal problem}. Multi-period portfolio selection problems are not considered here, however, one can expect a very similar behaviour to that observed in the single-period case. 
\par Finally, we could easily deal with alternative risk measures, such as the variance (e.g. \cite{JoFinance:Markowitz}), combination of CVaR and MAsD (e.g. see \cite{JFRM:BirungiMuthoni}), or approximations of other risk measures via multiple CVaR measures (see \cite{JBF:InuiKijima}). These were not included here for brevity of presentation.

\subsection{Penalized quantile regression}
\par Next we consider linear regression models of the following form
\[y_i = \beta_0 + \xi_i^\top \beta + \epsilon_i, \qquad i \in \{1,\ldots,l\}\]
\noindent where $\xi_i$ is a $d$-dimensional vector of covariates, $(\beta_0,\beta)$ are the regression coefficients and $\epsilon_i$ is some random error. A very popular problem in statistics is the estimation of the optimal coefficients, in the sense of minimizing a model of the following form:
\begin{equation} \label{eqn: linear regression model}
    \min_{(\beta_0,\beta)\ \in\ \mathbb{R} \times \mathbb{R}^d} \left\{\frac{1}{l} \sum_{i=1}^l \ell\left(y_i - \beta_0 - \xi_i^\top \beta\right) + \lambda p(\beta) \right\},
\end{equation}
\noindent where $\ell(\cdot)$ is some loss function and $p(\cdot)$ is a penalty function with an associated regularization parameter $\lambda \geq 0$. Following \cite{JRSS:Zou}, we consider the elastic-net penalty, 
\[p(\beta) \equiv p_{\tau}(\beta) \coloneqq \tau \|\beta\|_1 + \frac{1-\tau}{2}\|\beta\|_2^2,\qquad 0\leq \tau \leq 1.\]
\noindent For the loss function, we employ the quantile loss
\begin{equation} \label{eqn: quantile loss}
\ell(w) \equiv \rho_{\alpha}(w) \coloneqq \left(1-\alpha\right) w_- + \alpha w_+ = \frac{1}{2}\left(\lvert w\rvert + (2\alpha-1)w\right), \qquad 0 < \alpha < 1,
\end{equation}
\noindent where $w \in \mathbb{R}$. Notice that the case $\alpha = \frac{1}{2}$ yields the absolute loss. Letting $x = \begin{bmatrix} \beta_0 & \beta^\top\end{bmatrix}^\top$, and using Remark \ref{remark: l1 norm reformulation}, we can re-write problem \eqref{eqn: linear regression model} in the form of \eqref{primal problem}, as
\begin{equation*}
\begin{split}
\min_{(x,w)\ \in\ \mathbb{R}^{1+d}\times \mathbb{R}^l}&\ \left\{(\alpha-1)\mathds{1}_l^\top w + \frac{1}{2}x^\top Q x + \sum_{i=1}^l (w_i)_+  + \|Dx\|_1\right\},\\
\textnormal{s.t.}\qquad\ &\ \frac{1}{l}\left(-\begin{bmatrix}1 & \xi_i^\top \end{bmatrix}x + y_i \right) - w_i= 0,\qquad i = 1,\ldots,l,
\end{split}
\end{equation*}
\noindent where
\[Q = \begin{bmatrix}
 0 & 0_{1,d}\\ 0_{d,1} & \lambda (1-\tau) I_d
\end{bmatrix},\qquad D = \begin{bmatrix}
 0 & 0_{1,d}\\ 0_{d,1} & \lambda \tau I_d
\end{bmatrix}.\]
\paragraph{Real-world datasets:} In what follows, we solve several instances of problem  \eqref{eqn: linear regression model}. We showcase the effectiveness of the proposed approach on 5 regression problems taken from the LIBSVM library (see \cite{CC01a}). Additional information on the datasets is collected in Table \ref{Table: quantile regression datasets}. 
\begin{table}[!ht] 
\centering
\caption{Quantile regression datasets.\label{Table: quantile regression datasets}}
\scalebox{1}{
\begin{tabular}{lll}     
\specialrule{.2em}{.05em}{.05em} 
    \textbf{Name} & \textbf{\# of data points} & \textbf{\# of features} \\    \specialrule{.1em}{.3em}{.3em} 
space\_ga & 3107 & 6\\ 
abalone & 4177 & 8 \\ 
cpusmall & 8192 & 12 \\ 
cadata & 20,640 & 8 \\ 
E2006-tfidf & 16,087 &  150,360 \\ 
\specialrule{.2em}{.05em}{.05em} 
\end{tabular}}
\end{table}
\par We fix $\texttt{tol} = 10^{-4}$, $\lambda = 10^{-2}$, and $\tau = 0.5$, and run all three methods (active-set (AS), IP-PMM and OSQP) on the 5 instances for varying quantile level $\alpha$. The results are collected in Table \ref{Table quantile regression: varying confidence level}.
\begin{table}[!ht]
\caption{Quantile regression: varying confidence level (\texttt{tol} = $10^{-4}$, $\lambda = 10^{-2}$, $\tau = 0.5$).\label{Table quantile regression: varying confidence level}}
\centering
\begin{threeparttable}
\begin{tabular}{llllllll}     
\specialrule{.2em}{.05em}{.05em} 
    \multirow{2}{*}{\textbf{Problem}} & \multirow{2}{*}{$\bm{\alpha}$} & \multicolumn{3}{c}{\textbf{Iterations}}     & \multicolumn{3}{c}{\textbf{Time (s)}}  \\   \cmidrule(l{2pt}r{2pt}){3-5} \cmidrule(l{2pt}r{2pt}){6-8}
&   &  {{PMM}(SSN)[Fact.]} & {{IP--PMM}} & {{OSQ}} &  {{AS}} & {{IP--PMM}}    &{{OSQP}}  \\ \specialrule{.1em}{.3em}{.3em} 
      \multirow{4}{*}{{space\_ga}} & $0.50$ &  15(43)[30] & 19 & 2350 & 0.38 &  \textbf{0.35}  & 0.51\\ 
  & $0.65$ &  15(50)[35] & 20  & 16,300 & \textbf{0.37} &  0.39 & 3.36\\ 
 & $0.80$ &  15(62)[48] & 19 &  16,975 & 0.46 &  \textbf{0.35} & 3.46\\
   & $0.95$ &  15(78)[77] & 20  & 7400 & 0.48 &  \textbf{0.36}&  1.51 \\ \specialrule{.00002em}{.1em}{.1em}
 \multirow{4}{*}{{abalone}} & $0.50$ &  14(68)[64] & 35 & 6350 & \textbf{1.02} &  1.29 &  2.15\\ 
  & $0.65$ &  14(68)[63] & 35 & 6300 & \textbf{0.89} &  1.36 & 2.17\\ 
 & $0.80$ &  14(79)[74] & 27 & 12,425  & \textbf{0.89} &  1.01& 4.13 \\
   & $0.95$ &  14(87)[79] & 18 & 2875& 0.82 &  \textbf{0.69} & 0.97 \\ \specialrule{.00002em}{.1em}{.1em}
    \multirow{4}{*}{{cpusmall}} & $0.50$ &  15(64)[64] & 29 & $\doublebarwedge$\tnote{1} & \textbf{2.01} &  2.57 & $\doublebarwedge$ \\ 
  & $0.65$ &  16(74)[74] & 30 & $\doublebarwedge$ & \textbf{2.33} &  2.81 & $\doublebarwedge$\\ 
 & $0.80$ &  16(86)[85] & 26 & $\doublebarwedge$ & 2.64 &  \textbf{2.37}&  $\doublebarwedge$\\
   & $0.95$ &  15(110)[110] & 22 & $\doublebarwedge$ & 2.85 &  \textbf{1.96} & $\doublebarwedge$\\ \specialrule{.00002em}{.1em}{.1em}
 \multirow{4}{*}{{cadata}} & $0.50$ &  3(66)[65] & 49 &  7125 & \textbf{2.96} &  14.63 & 16.10 \\ 
  & $0.65$ &  4(85)[83] & 42 & 9050  & \textbf{3.65} &  12.43 & 20.27  \\ 
 & $0.80$ &  3(77)[76] & 45 & 7275 & \textbf{3.45} &  13.13 &   16.45 \\
   & $0.95$ &  3(225)[224] & 80 & 7450 & \textbf{13.14} &  23.32 & 17.08 \\ \specialrule{.00002em}{.1em}{.1em}
\multirow{4}{*}{{E2006-tfidf}} & $0.50$ &  14(27)[20] & $\dagger$\tnote{2} & $\dagger$  & \textbf{84.82} &  $\dagger$  & $\dagger$\\ 
  & $0.65$ &  15(34)[26] & $\dagger$  & $\dagger$ & \textbf{95.08} & $\dagger$ & $\dagger$\\ 
 & $0.80$ &  16(46)[35] & $\dagger$ & $\dagger$ & \textbf{112.24} &  $\dagger$ & $\dagger$ \\
   & $0.95$ &  17(79)[71] & $\dagger$ & $\dagger$ & \textbf{165.78} &  $\dagger$ &  $\dagger$\\ 
\specialrule{.2em}{.05em}{.05em} 
\end{tabular}
\begin{tablenotes}
\item[1] $\doublebarwedge$ indicates that the solver incorrectly identified the problem as infeasible.
\item[2] $\dagger$ indicates that the solver ran out of memory.
\end{tablenotes}
\end{threeparttable}
\end{table}
\par From Table \ref{Table quantile regression: varying confidence level} we observe that the three approaches are comparable for the smaller isntances (space\_ga, abalone, and cpusmall). The active-set scheme significantly outperforms IP-PMM and OSQP on the cadata problem, mainly due to its better numerical stability. Indeed, this problem is highly ill-conditioned, and this is reflected in the increased number of IP-PMM iterations needed to obtain a 4-digit accurate solution. Finally, we observe that for the large-scale instance (E2006-tfidf), IP-PMM and OSQP crashed due to memory requirements. In this case, the active set scheme was warm-started by a matrix-free ADMM (as discussed in Section \ref{sec: warm start}), since the standard ADMM also crashed due to excessive memory requirements. The active-set scheme, however, manages to solve these large-scale instances very efficiently, without running into any memory issues (consistently). The OSQP solver is competitive for certain instances, albeit slightly slower than both second-order solvers, but fails to solve the cpusmall instances due to an incorrect classification of these instances as infeasible. Additionally, we observe that it is less consistent than the second-order approaches, since its iterations vary greatly with the problem parameters (see instances space\_ga and abalone). 
\par Next, we fix $\texttt{tol} = 10^{-4},$ and $\alpha = 0.8$, and run the three methods for varying regularization parameters $\tau$ and $\lambda$. The results are collected in Table \ref{Table quantile regression: varying regularization}. We observe that both second-order methods are quite robust for a wide range of parameter values. OSQP is competitive for certain instances, however, its behaviour is greatly influenced by the problem parameters. Nonetheless, it was albe to outperform the second-order solvers on two out of the four cadata instances. Finally, the active-set scheme consistently solved the large-scale instance (E2006-tfidf) for a wide range of parameter values, without running into any memory issues. 

\begin{table}[!ht]
\caption{Quantile regression: varying regularization (\texttt{tol} = $10^{-4}$, $\alpha = 0.8$).\label{Table quantile regression: varying regularization}}
\centering
\scalebox{0.94}{\begin{threeparttable}
\begin{tabular}{llllllll}     
\specialrule{.2em}{.05em}{.05em} 
    \multirow{2}{*}{\textbf{Problem}} & \multirow{2}{*}{$(\bm{\tau},\bm{\lambda})$} & \multicolumn{3}{c}{\textbf{Iterations}}     & \multicolumn{3}{c}{\textbf{Time (s)}}  \\   \cmidrule(l{2pt}r{2pt}){3-5} \cmidrule(l{2pt}r{2pt}){6-8}
&   &  {{PMM}(SSN)[Fact.]} & {{IP--PMM}} & {{OSQP}}  &  {{AS}} & {{IP--PMM}}   &{{OSQP}}   \\ \specialrule{.1em}{.3em}{.3em} 
      \multirow{4}{*}{{space\_ga}}& $(0.2,5\cdot 10^{-2})$ &  15(58)[49] & 16 & 37,125 & 0.51 &  \textbf{0.30} & 8.00\\
      &$(0.4,1\cdot 10^{-2})$ &  15(61)[48] & 26 & 12,050 & \textbf{0.46} &  0.49 & 2.48\\
       &$(0.6,5\cdot 10^{-3})$ &  15(67)[46] & 25 & 13,025 & \textbf{0.44} &  0.46 & 2.68 \\
        &$(0.8,1\cdot 10^{-3})$ &  15(82)[62] & 27 & 20,575 & 0.61 &  \textbf{0.49} & 4.21\\ \specialrule{.00002em}{.1em}{.1em}
 \multirow{4}{*}{{abalone}} & $(0.2,5\cdot 10^{-2})$ &  14(70)[65] & 25  & 16,225 & \textbf{0.88} &  1.32 & 5.49 \\
      &$(0.4,1\cdot 10^{-2})$ &  14(71)[67] & 34  & 15,575 & \textbf{0.86} &  1.48 & 5.16\\
       &$(0.6,5\cdot 10^{-3})$ &  15(80)[76] & 39 & 18,675 & \textbf{0.92} &  1.50 & 6.26\\
        &$(0.8,1\cdot 10^{-3})$ &  15(116)[109] & 30 & 9975 & 1.44 &  \textbf{1.16} & 3.34 \\ \specialrule{.00002em}{.1em}{.1em} 
 \multirow{4}{*}{{cpusmall}} & $(0.2,5\cdot 10^{-2})$ &  15(83)[81] & 20 & $\doublebarwedge$\tnote{1} & \textbf{1.98} &  2.14 & $\doublebarwedge$ \\
      &$(0.4,1\cdot 10^{-2})$ &  15(81)[78] & 26& $\doublebarwedge$  & \textbf{2.10} &  2.39 & $\doublebarwedge$\\
       &$(0.6,5\cdot 10^{-3})$ &  16(101)[94] & 21 & $\doublebarwedge$  & 2.79 &  \textbf{1.94} & $\doublebarwedge$ \\
        &$(0.8,1\cdot 10^{-3})$ &  15(106)[103] & 20 & $\doublebarwedge$ & 3.16 &  \textbf{1.85} & $\doublebarwedge$\\ \specialrule{.00002em}{.1em}{.1em} 
 \multirow{4}{*}{{cadata}} & $(0.2,5\cdot 10^{-2})$ &  9(147)[143] & 62 &  5875 & 15.42 &  18.05 & \textbf{13.18}  \\
      &$(0.4,1\cdot 10^{-2})$ &  4(63)[62] & 50 & 7350 & \textbf{2.76} &  14.68 & 17.41\\
       &$(0.6,5\cdot 10^{-3})$ &  5(121)[118] & 45 & 9775 & \textbf{5.49} &  13.32 & 21.70\\
        &$(0.8,1\cdot 10^{-3})$ &  7(276)[276] & 33 & 875 & 14.01 &  9.91 & \textbf{2.08} \\ \specialrule{.00002em}{.1em}{.1em}
\multirow{4}{*}{{E2006-tfidf}} & $(0.2,5\cdot 10^{-2})$ &  16(47)[38] & $\dagger$\tnote{2} & $\dagger$ & \textbf{115.87} &  $\dagger$ & $\dagger$\\ 
  & $(0.4,1\cdot 10^{-2})$ &  16(43)[36] & $\dagger$ & $\dagger$&  \textbf{118.89} & $\dagger$ & $\dagger$ \\ 
 & $(0.6,5\cdot 10^{-3})$ &  16(43)[34] & $\dagger$ & $\dagger$ & \textbf{115.90} &  $\dagger$ & $\dagger$ \\
   & $(0.8,1\cdot 10^{-3})$ &  17(48)[40] & $\dagger$ & $\dagger$ & \textbf{128.38} &  $\dagger$ & $\dagger$\\
\specialrule{.2em}{.05em}{.05em} 
\end{tabular}
\begin{tablenotes}
\item[1] $\doublebarwedge$ indicates that the solver incorrectly identified the problem as infeasible.
\item[2] $\dagger$ indicates that the solver ran out of memory.
\end{tablenotes}
\end{threeparttable}}
\end{table}

\subsection{Binary classification via linear support vector machines}
\par Finally, we are interested in training a binary linear classifier using regularized soft-margin linear support vector machines (SVMs), \cite{Springer:Vapnik}. More specifically, given a training dataset  $\{(y_i,\xi_i)\}_{i=1}^l$, where $y_i \in \{-1,1\}$ are the \emph{labels} and $\xi_i \in \mathbb{R}^d$ are the \emph{feature vectors} (with $d$ the number of features), we would like to solve the following optimization problem
\begin{equation} \label{eqn: linear SVM}
\min_{(\beta_0,\beta) \in \mathbb{R}\times\mathbb{R}^d} \left\{\frac{1}{l}\sum_{i=1}^l \left( 1- y_i\left(\xi_i^\top \beta -\beta_0\right)\right)_+  + \lambda \left(\tau_1 \|\beta\|_1 + \frac{\tau_2}{2}\|\beta\|_2^2\right)\right\},
\end{equation}
\noindent where $\lambda > 0$ is a regularization parameter, and $\tau_1,\ \tau_2 > 0$ are the weights for the $\ell_1$ and $\ell_2$ regularizers, respectively. The standard Euclidean regularization is traditionally used as a trade-off between the margin of the classifier (the larger the better) and the correct classification of $\xi_i$, for all $i \in \{1,\ldots,l\}$ (e.g. \cite{MachLear:CortesVapnik}). However, this often leads to a dense estimate for $\beta$, which has led researchers in the machine learning community into considering the $\ell_1$ regularizer instead, in order to encourage sparsity (e.g. \cite{ICML:BradManga}). It is well-known that both regularizers can be combined to obtain the effect of both individual ones, using the elastic-net penalty (see, for example, \cite{StatSin:Wang_etal}), assuming that $\tau_1$ and $\tau_2$ are appropriately tuned.

\paragraph{Real-world datasets:} In what follows we consider elastic-net SVM instances of the form of \eqref{eqn: linear SVM}. We showcase the effectiveness of the proposed active-set scheme on 3 large-scale binary classification datasets taken from the LIBSVM library (\cite{CC01a}). The problem names, as well as the number of features and training points are collected in Table \ref{Table: Binary classification datasets}.
\begin{table}[!ht] 
\centering
\caption{Binary classification datasets.\label{Table: Binary classification datasets}}
\scalebox{1}{
\begin{tabular}{lll}     
\specialrule{.2em}{.05em}{.05em} 
    \textbf{Name} & \textbf{\# of training points} & \textbf{\# of features} \\  \specialrule{.1em}{.3em}{.3em} 
rcv1 & 20,242 & 47,236 \\ 
real-sim & 72,309 & 20,958\\
news20 & 19,996 &  1,355,191\\
\specialrule{.2em}{.05em}{.05em} 
\end{tabular}}
\end{table}
\par In the experiments to follow, we only consider large-scale problems that neither IP-PMM nor OSQP can solve, due to excessive memory requirements. Nonetheless, we should note that, by following the developments in \cite{SIREV:DeSimone_etal}, IP-PMM can be specialized to problems of this form in order to be able to tackle large-scale instances. However, the same can be said about the proposed active-set method. The schemes tested in this work employ factorization for the solution of their associated linear systems. The introduction of preconditioning, e.g. as in \cite{arxiv:PougkGondzio}, could improve their efficiency in certain cases, but is out of the scope of this article. In the following experiments, we employ the (matrix-free) prox-linear ADMM (as described in Section \ref{sec: warm start}), to avoid any memory issues in the warm-starting phase of the proposed algorithm.
\par We fix $\texttt{tol} = 10^{-5}$, $\lambda = 10^{-2}$, and run the active-set solver for all datasets given in Table \ref{Table: Binary classification datasets} for varying values of the regularization parameters $\tau_1$, and $\tau_2$. The results are collected in Table \ref{Binary classification via elastic-SVMs: varying regularization }.
\begin{table}[!ht]
\caption{Binary classification via elastic-net SVMs: varying regularization (\texttt{tol} = $10^{-5}$, $\lambda = 10^{-2}$).\label{Binary classification via elastic-SVMs: varying regularization }}
\centering
\scalebox{1}{\begin{threeparttable}
\begin{tabular}{lllll}     
\specialrule{.2em}{.05em}{.05em} 
    \multirow{2}{*}{\textbf{Problem}} & \multirow{2}{*}{$\bm{\tau_1}$} & \multirow{2}{*}{$\bm{\tau_2}$} & \multicolumn{1}{c}{\textbf{Iterations}}     & \multicolumn{1}{c}{\textbf{Time (s)}}  \\   \cmidrule(l{2pt}r{2pt}){4-4} \cmidrule(l{2pt}r{2pt}){5-5}
&  &  &  {{PMM}(SSN)[Fact.]}    &  {{AS}}    \\ \specialrule{.1em}{.3em}{.3em} 
     \multirow{4}{*}{{rcv1}}&  0.2 & 0.2 &  48(139)[105] & 43.36\\
      & 0.8 & 0.2 &  47(100)[58] & 23.28 \\
          &  0.2 & 0.8 &  47(165)[133] & 54.81 \\
      & 5 & 5 &  39(85)[45] & 20.57 \\ \specialrule{.00002em}{.1em}{.1em}
              \multirow{4}{*}{{real-sim}}&  0.2 & 0.2 &  50(130)[100] & 252.25 \\
      &  0.8 & 0.2 &  47(85)[48] & 132.20 \\
     &  0.2 & 0.8  &  48(101)[68] & 201.56\\
      & 5 & 5 &  47(85)[48] & 128.75 \\  \specialrule{.00002em}{.1em}{.1em}
               \multirow{4}{*}{{news20}} &  0.2 & 0.2 &  50(142)[102] & 213.12 \\
      &  0.8 & 0.2 &  41(85)[43] & 133.94 \\
         &  0.2 & 0.8  &  70(212)[161] & 74.42 \\
      & 5 & 5 &  41(85)[43] & 130.33 \\
\specialrule{.2em}{.05em}{.05em} 
\end{tabular}
\end{threeparttable}}
\end{table}
\par We observe that the solver was able to consistently solve these instances, without running into memory issues. In this case, the behaviour of the active-set solver was affected by the parameters $\tau_1$ and $\tau_2$ (indeed, see the number of SSN iterations for different regularization values), however, we consistently obtain convergence in a very reasonable amount of time. Overall, we observe that the proposed algorithm is very general and can be applied in a plethora of very important applications arising in practice. We were able to showcase that the active-set nature of the method allows one to solve large-scale instances on a personal computer, without the need of employing iterative linear algebra (which could also complement the solver, as in \cite{arxiv:PougkGondzio}). The proposed method strikes a good balance between first-order methods, which are fast but unreliable, and second-order interior point methods, which are extremely robust but can struggle with the problem size, ill-conditioning and memory requirements. We have demonstrated that $\ell_1$-regularized convex quadratic problems with piecewise-linear terms can be solved very efficiently using the proposed active-set scheme, and we conjecture that the algorithm can be readily extended to deal with general nonlinear convex objectives, or discretizations of stochastic two-stage problems. 

\iffalse
\subsection{Alternative applications}
\par The problem in \eqref{primal problem} can readily model a plethora of very important applications arising in practice. Immediate alternative applications for which the proposed solver is expected to be especially efficient include, among many others, classification via the hinge regression (e.g. see \cite{NeurComp:Rosasco}), discretized risk-averse and/or $L^1$-regularized partial differential equation optimization (e.g. see \cite{SIAMOpt:Kouri,arxiv:PougkGondzio}), holding and backlogging cost in inventory management (e.g. see . Additionally the method can be used as an inner solver for nonlinear CVaR or MAsD minimization problems solved via the Gauss-Newton method \cite{MathProg:Burke}. Finally, it can be used as an inner solver for sub-problems arising within combinatorial optimization solvers involving relaxed models of the form of \eqref{primal problem}.
\fi

\section{Conclusions} \label{sec: Conclusions}
\par In this paper we derived an efficient active-set method for the solution of convex quadratic optimization problems with piecewise-linear terms in the objective. The method, which complements our developments in the accompanying paper [``\emph{An active-set method for sparse approximations. Part I: Separable $\ell_1$ terms}", \emph{S. Pougkakiotis,  J. Gondzio, D. S. Kalogerias}], arises by suitably combining a proximal method of multipliers with a semismooth Newton scheme, and admits an active-set interpretation. By taking advantage of the piecewise-linear terms in the objective, the method has very reasonable memory requirements since it utilizes only a small active-set at every inner-outer iteration. We warm-start the algorithm using an appropriate alternating direction method of multipliers, and ensure faster convergence and reduced memory requirements. 
\par We showcase the efficiency and robustness of the proposed scheme on a variety of optimization problems arising in risk-averse portfolio optimization, quantile regression, and binary classification via linear support vector machines. A numerical comparison against a robust interior point method and a state-of-the-art alternating direction method of multipliers demonstrates the viability of the approach as well as its limited memory requirements. In particular, we observe a significantly better behaviour, compared to the two other solvers, when dealing with large-scale instances. Overall, the approach remains efficient for a wide range of problems and strikes a good balance between cheap but unreliable first-order methods and expensive but highly reliable interior point methods.

\appendix
\section{Appendix: Termination criteria} \label{Appendix: termination criteria}
\par The optimality conditions of \eqref{primal problem} can be written as
\begin{equation*} 
\begin{split}
\textbf{prox}_{g_1}\left(x -c - Qx + \begin{bmatrix}C^\top & A^\top \end{bmatrix} y - z\right) =\ x,&\qquad
\textbf{prox}_{g_2}\left(w - y_{1:l}\right) =\ w,\\
\begin{bmatrix} Cx + d -  w\\ Ax-b \end{bmatrix} = \ 0_{l+m},&\qquad
 \Pi_{\mathcal{K}}(x + z) =\ x,
\end{split}
\end{equation*}
\noindent and the termination criteria for Algorithm \ref{PMM algorithm} (given a tolerance $\epsilon > 0$) can be summarized as
\begin{equation} \label{termination criteria for PD-PMM}
\begin{split}
\frac{\left\|x - \textbf{prox}_{g_1}\left(x -c - Qx + \begin{bmatrix}C^\top & A^\top \end{bmatrix} y - z\right)\right\|}{1+\|c\|_{\infty}} \leq \epsilon,&\quad
\left\|w - \textbf{prox}_{g_2}\left(w - y_{1:l}\right)\right\| \leq \epsilon,\\
\frac{\left\|\begin{bmatrix} Cx + d -  w\\ Ax-b \end{bmatrix}\right\|}{1+\|b\|_{\infty} + \|d\|_{\infty}} \leq \epsilon,&\quad
\frac{\|x - \Pi_{\mathcal{K}}(x + z)\|}{1+\|x\|_{\infty} + \|z\|_{\infty}} \leq \epsilon.
\end{split}
\end{equation}
\par From the reformulation of \eqref{primal problem} given in \eqref{primal problem ADMM reformulation}, the termination criteria of Algorithm \ref{proximal ADMM algorithm} are as follows (upon noting that the variables of the algorithm are $(x,w,u,y)$)
\begin{equation} \label{termination criteria for pADMM}
\begin{split}
\frac{\left\|c + Qx - \begin{bmatrix} C^\top & A^\top -I_n & 0_{n,l}\end{bmatrix} y\right\|}{1+\|c\|} \leq \epsilon,&\qquad
\left\|\begin{bmatrix} I_l & 0_{l,m+n} & I_l\end{bmatrix}y\right\| \leq \epsilon,\\
\frac{\left\|M_r \begin{bmatrix}
x\\
w\\
u
\end{bmatrix} - \begin{bmatrix}
-d\\
b\\
0_{l+n}
\end{bmatrix}\right\|}{\left\|\begin{bmatrix}
-d\\ b\end{bmatrix}\right\|+1}
 \leq \epsilon,&\qquad 
\frac{\left\|u - \Pi_{\mathcal{K}\times \mathbb{R}^l}\left(\textbf{prox}_g\left( u +\tilde{y}\right) \right) \right\|}{1+ \|u\|+\left\|\tilde{y}\right\|} \leq \epsilon,
\end{split}
\end{equation}
\noindent where $\tilde{y} \coloneqq y_{(l+m+1:2l+m+n)}$.
\bibliography{references} 
\bibliographystyle{siamplain}

\end{document}